\documentclass[aos,secthm,seceqn,dvips]{arximspdf}
\usepackage{graphicx}

\doi{10.1214/009053607000000046}
\volume{35}
\issue{4}
\pubyear{2007}
\firstpage{1432}
\lastpage{1455}

\makeatletter
\newcommand{\ed}{\iffalse}

\newtheorem{lem}{Lemma}[section]
\newproclaim{rem}{Remark}[section]
\newproclaim{example}{Example}[section]
\makeatother

\begin{document}
\begin{frontmatter}

\title{Dependency and false discovery rate: Asymptotics\thanksref{AA}}
\pdftitle{Dependency and false discovery rate: Asymptotics}
\thankstext{AA}{Supported by the Deutsche Forschungsgemeinschaft.}
\runtitle{Dependency and FDR: Asymptotics}

\begin{aug}
\author[A]{\fnms{Helmut} \snm{Finner}\corref{}\ead[label=e1]{finner@ddz.uni-duesseldorf.de}},
\author[A]{\fnms{Thorsten} \snm{Dickhaus}\ead[label=e2]{dickhaus@ddz.uni-duesseldorf.de}}
and
\author[B]{\fnms{Markus} \snm{Roters}\ead[label=e3]{Markus.Roters@omnicarecr.com}}
\runauthor{H. Finner, T. Dickhaus and M. Roters}
\affiliation{German Diabetes Center
D\"usseldorf, German Diabetes Center D\"{u}sseldorf\break and Omnicare Clinical Research}
\address[A]{
H. Finner\\
T. Dickhaus\\
German Diabetes Center\\
\quad at the Heinrich--Heine--Universit\"at
D\"usseldorf\\
Institute of Biometrics and Epidemiology\\
D\"usseldorf\\
Germany\\
\printead{e1}\\
\phantom{E-mail: }\printead*{e2}} %adresu isvedimo komanda gale!
\address[B]{
M. Roters\\
Omnicare Clinical Research\\
Biometrics Department\\
K\"oln\\
Germany\\
\printead{e3}}
\end{aug}

% HISTORY:
\received{\smonth{4} \syear{2006}}
\revised{\smonth{9} \syear{2006}}

% ABSTRACT
%
\begin{abstract}
Some effort has been undertaken over the last decade to provide
conditions for the control of the false discovery rate by the linear
step-up procedure (LSU) for testing $n$ hypotheses when test statistics
are dependent. In this paper we investigate the expected error rate
(EER) and the false discovery rate (FDR) in some extreme parameter
configurations when $n$ tends to infinity for test statistics being
exchangeable under null hypotheses. All results are derived in terms of
$p$-values. In a general setup we present a series of results concerning
the interrelation of Simes' rejection curve and the (limiting)
empirical distribution function of the $p$-values. Main objects under
investigation are largest (limiting) crossing points between these
functions, which play a key role in deriving explicit formulas for EER
and FDR. As specific examples we investigate equi-correlated normal and
$t$-variables in more detail and compute the limiting EER and FDR
theoretically and numerically. A surprising limit behavior occurs if
these models tend to independence.
\end{abstract}

% KEYWORDS
%
\begin{keyword}[class=AMS]
\kwd[Primary ]{62J15}
\kwd{62F05}
\kwd[; secondary ]{62F03}
\kwd{60F99}.
\end{keyword}
\begin{keyword}
\kwd{Exchangeable test statistics} \kwd{expected error rate} \kwd{false discovery rate}
\kwd{Glivenko--Cantelli theorem} \kwd{largest crossing point} \kwd{least favorable
configurations} \kwd{multiple comparisons} \kwd{multiple test
procedure}
\kwd{multivariate total positivity of order 2} \kwd{positive regression
dependency} \kwd{Simes' test}.
\end{keyword}

\end{frontmatter}

%s1 ###
\section{Introduction}

Control of the false discovery rate (FDR) in multiple hypotheses
testing has become an attractive approach especially if a large number
of hypotheses is at hand. The first FDR controlling procedure, a linear
step-up procedure (LSU), was originally designed for independent
$p$-values (cf. \cite{r1}) and has its origins in \cite{r3} (cf. also
\cite{r15}). Meanwhile, it is known that the LSU-procedure controls the
FDR even if the test statistics obey some special dependence structure.
Key words are MTP$_2$ (multivariate total positivity of order 2) and
PRDS (positive regression dependency on subsets).
More formal descriptions of these conditions and proofs can be found in
\cite{r2}~and~\cite{r14}. In view of testing problems with some ten
thousand hypotheses as they appear, for example, in genetics,
asymptotic considerations become more and more popular. The first
asymptotic investigations concerning expected type I errors of the
LSU-procedure, as well as for the corresponding linear step-down (LSD)
procedure for the independent case, can be found in
\cite{r9}~and~\cite{r10}. A first theoretical comparison of classical
stepwise procedures controlling a multiple level $ \alpha$ [or
familywise error rate (FWER) in the strong sense] based on asymptotics
is given in~\cite{r5}.
Moreover, attempts to improve the LSU-procedure and interesting
investigations based on asymptotics can be found, for example, in
\mbox{\cite{r12,r11}} and \cite{r22}.

The LSU-procedure is based on the critical values $\alpha_{i:n} = i
\alpha/ n $, $ i = 1 , \ldots, n $, introduced in \cite{r17} in a
different context. Based on ordered $p$-values\break $ p_{1:n} \leq \cdots
\leq p_{n:n} $, the LSU-procedure rejects the corresponding
hypotheses $ H_{1:n} , \ldots, H_{m:n} $, where $ m = \max\{i\dvtx
p_{i:n} \leq\alpha_{i:n} \} $. The corresponding LSD-pro\-cedure
rejects $ H_{1:n} , \ldots, H_{r:n} $, where $ r = \max\{ i\dvtx
p_{j:n} \leq\alpha_{j:n} \mbox{ for all } j = 1 , \ldots, i \} $.
Since $ m \geq r $, the LSU-procedure may reject more hypotheses than
the LSD-procedure, never less. In this paper we restrict attention to
the LSU-procedure, which can be rewritten in terms of the empirical
c.d.f. (e.c.d.f.) $ {F}_n $ (say) of the~$p_i$'s. Setting $ t^* = \sup
\{ t\dvtx{F}_n ( t ) \geq t / \alpha\} $, $ H_i $ is rejected iff $
p_i \leq t^* $. The rejection curve $ r_\alpha( t ) = t / \alpha$
will be called the Simes-line. Note that $ \alpha_{i: n} =
r_\alpha^{-1} ( i/n ) $. The threshold~$ t^* $ will be called the
largest crossing point (LCP) between the e.c.d.f. and the Simes-line and plays
a crucial role in this paper.

FDR control for a multiple test procedure is defined as follows. Let $
V_n $ denote the number of falsely rejected null hypotheses and let $
R_n $ denote the number of all rejections. Then the FDR (depending on
the underlying parameter configuration $ \vartheta\in\Theta$, say) is
defined by
\[
\operatorname{FDR}_n ( \vartheta) = \mathbb{E}\biggl[ \frac
{V_n}{R_n \vee1 }\biggr]
\]
and is said to be controlled at level $ \alpha$ if
\[
\sup_{\vartheta\in\Theta} \operatorname{FDR}_n ( \vartheta) \leq
\alpha.
\]
The ratio $ V_n / [ R_n \vee1 ] $ is the false discovery proportion
(FDP).
In the case of independent $p$-values both LSU and LSD control the FDR at
level $ \alpha$; more precisely, if $ \vartheta\in\Theta$ is such that
exactly $ n_0 $ hypotheses are true and the remaining $n_1=n-n_0$ ones
are false, for both LSU and LSD, the actual FDR is bounded by $ n_0
\alpha/ n $. Under weak additional assumptions, we have in this setting
for the LSU-procedure %
\[
\operatorname{FDR}_n ( \vartheta) = \frac{n_0}{n} \alpha.
\]
Different proofs for this fact can be found in
\cite{r1,r9,r14} and \cite{r22}.

In \cite{r10} the expected number of type I errors (ENE), that is,
$\mathrm{ENE}_n (\vartheta) = \mathbb{E}[ V_n ] $, of LSU and LSD was
investigated for the case that all hypotheses are true and \mbox{$p$-values} are
independent. In this case the limiting ENE for $ n \to\infty$ equals $
\alpha/ ( 1 - \alpha)^2 $ for LSU and $ \alpha/ ( 1 - \alpha) $ for
LSD. Moreover, in \cite{r9} the expected type I error rate (EER)
defined by $ \operatorname{EER}_n ( \zeta) = \mathbb{E}[ V_n / n] $ was
studied if a proportion $ 1 - \zeta$ of hypotheses is {\em totally
false}, that is, with $p$-values equal to zero with probability $1$. For
independent $p$-values, it was shown in \cite{r9} under quite general
assumptions that for both LSU and LSD
\[
\lim_{ n \to\infty} \sup_{\vartheta\in\Theta} \mathbb{E}\biggl[
\frac{V_n}{n} \biggr]= \bigl(1- \sqrt{1-\alpha}\bigr)^2 / \alpha= \alpha/4 +
\alpha^2/8 + O(\alpha^3) \approx\alpha/4 .
\]
The worst case for the EER appears if the proportion of true hypotheses
tends to $ \zeta= ( 1 - \sqrt{ 1 - \alpha} )/ \alpha = 1 / 2 + \alpha/
8 + O(\alpha^2) $, and, for small values of $\alpha$, the expected type
I error rate is then approximately $ \alpha/ 4 $.

In this paper we investigate the behavior of EER and FDR of the
LSU-procedure based on dependent test statistics if the number of
hypotheses tends to infinity. It will be assumed that test statistics
are exchangeable under the corresponding null hypotheses. The main
issue will be the calculation of the limiting values of the actual EER
and FDR in some extreme parameter configurations, where a proportion
$\zeta_n$ of hypotheses will be assumed to be true and the remaining
hypotheses will be assumed to be totally false. These configurations
are least favorable for the EER, that is, EER becomes largest under
these configurations if $ \zeta_n$ is given.
Theoretical results on least/most favorable configurations for the FDR
(configurations where the FDR becomes largest/smallest) under
dependence remain a challenging open problem. However, simulations
indicate that extreme configurations ($n_0$ hypotheses true, $n_1$
hypotheses totally false) are first candidates for least favorable
configurations and therefore of special theoretical interest.
Until now, not many results are available concerning the behavior of
EER and FDR under dependence. A brief discussion on expected type I
errors for single-step procedures based on exchangeable test statistics
and range statistics can be found in~\cite{r8}.

In Section \ref{sec3} we develop a general theory for the computation
of the limiting EER and FDR assuming that exchangeable test statistics
of the type $ T_i = g ( X_i , Z ) $ are at hand. The results heavily
depend on the limit behavior of the e.c.d.f. of the underlying $p$-values
given the value $z$ of the disturbance variable $Z$. Generally, the
limiting e.c.d.f. $ F_\infty$ (say) of dependent $p$-values differs
substantially from that of independent $p$-values. Formulas for the
limiting e.c.d.f. and crossing point determination are summarized in
Lemma~\ref{lem3}. For $ \zeta_n \to \zeta \not= 1 $, limiting
$\operatorname{EER}$ and $\operatorname{FDR}$ are computed in Theorems
\ref{thm3_1}~and~\ref{thm3_2} in terms of the set of largest crossing
points (LCP's) between $ F_\infty$ and the Simes-line. The case $
\zeta_n \to 1 $ is more complex because limiting LCP's can be zero. For
the latter case, we derive some important technical results for the FDR
in Lemmas \ref{FDR_lem1}~and~\ref{lem3_4} supposing that the c.d.f. of a
proportion of $p$-values is linear in a neighborhood of zero.
%Results for the limiting behavior of $F_\infty$ in zero are given in
%Lemmas
The limiting $\operatorname{EER}$ and $\operatorname{FDR}$ are then
computed in Theorems \ref{thm2_3}~and~\ref{thm2_4}. Moreover, we give an
example where the FDR is exactly the same as in the independent case.
By utilizing the results of Section \ref{sec3}, we investigate
equi-correlated normal variables in Section \ref{sec4} and jointly
studentized $t$-statistics in Section \ref{sec5}. The corresponding
formulas for the limiting $\operatorname{EER}$ and $\operatorname{FDR}$
are given in Theorems \ref{thm4_1} and \ref{thm4_2} and Theorems
\ref{thm5_1}~and~\ref{thm5_2}, respectively. A surprising behavior of the
FDR occurs if these models tend to independence and the proportion of
false hypotheses tends to~$0$; see Theorem \ref{thm4_3} and Theorem
\ref{F_nu}. Some figures in Sections \ref{sec4} and \ref{sec5}
illustrate the limiting behavior of EER and FDR. The numerical and
computational effort for these graphs was enormous. A few concluding
remarks are given in Section \ref{sec6}. Short proofs are in the main
text, while more technical proofs are deferred to the
\hyperref[app]{Appendix}.

%s2 ###
\section{Exchangeable test statistics: general considerations} \label{sec3}

We first consider the following basic model with exchangeable test
statistics. Let $ X_i, \; i=1, \ldots, n $, be real-valued independent
random variables with support $ \mathcal{X}$. Moreover, let $ Z $ be a
further real-valued random variable, independent of the $ X_i$'s, with
support $ \mathcal{Z}$ and continuous c.d.f. $ W_Z $. Denote the c.d.f.
of $ X_i $ by $ W_i $.
Suppose the c.d.f. $ W_i $ depends on a parameter $ \vartheta_i \in[
\vartheta_0 , \infty) $, where $ \vartheta_0 $ is known. Without loss
of generality, it will be assumed that $ \vartheta_0 = 0 $. Consider the
multiple testing problem
\[
H_i\dvtx \vartheta_i = 0 \quad \mbox{versus} \quad K_i \dvtx
\vartheta_i > 0,\qquad  i = 1 , \ldots, n.
\]
Suppose that $ T_i = g ( X_i , Z ) $ (with support $ \mathcal{T}$) is a
suitable real-valued test statistic for testing $ H_i $ such that $ T_i
$ tends to larger values if $ \vartheta_i $ increases. In
Section~\ref{sec4} we consider statistics of the type $ T_i = g ( X_i ,
Z ) = X_i - Z $ and in Section~\ref{sec5} $ T_i = g ( X_i , Z ) = X_i /
Z $; see Examples \ref{Example_2.1}~and~\ref{Example_2.2} below. The
sets $\mathcal{X}$, $\mathcal{Z}$ and $\mathcal{T}$ are assumed to be
intervals. For convenience, we assume in this section that $ g $ is
continuous, strictly increasing in the first and either strictly
monotone or constant in the second argument. Moreover, let $ g_1 $
denote the inverse of $ g $ with respect to the first argument of $ g
$, that is, $ g ( x , z ) = w \mbox{ iff } x = g_1 ( w , z ) $. If $g$
is strictly monotone in the second argument, we denote the inverse of $
g $ with respect to the second argument by $g_2$, that is, $ g ( x , z
) = w \mbox{ iff } z = g_2 ( x , w ) $.

In the case that $ H_i $ is true, the c.d.f. of $ X_i $ ($ T_i $) will
be denoted by $ W_X $ ($ W_T $) and $ W_X $ is assumed to be
continuous. For $ Z = z $, we define $p$-values $ p_i = p_i (z) $ as a
function of $ z$ by
%e1 ###
%
\begin{equation} \label{p_value}
p_i ( z ) = 1 - W_T ( g ( x_i , z ) ) ,\qquad  i = 1 , \ldots, n.
\end{equation}
The ordered $p$-values are given by $p_{i:n} (z) = 1 - W_T ( g (
x_{n-i+1:n} , z ) ) $. Under\break $H_0=\bigcap_{i=1}^n H_i $, the e.c.d.f.
of the $p$-values is denoted by  $ F_n ( \cdot| z ) $.

\begin{rem} \label{rem2_1}
It is important to note that, given $ Z=z $, the $p$-values $ p_i (z) $,
$ i= 1 , \ldots, n $, can be regarded (a) as conditionally independent
random variables $ 1 - W_T ( g ( X_i , z ) ) $ with values in $[0,1]$,
or, (b) under $H_0$, as realizations of conditionally i.i.d. random
variables with a common c.d.f. $ F_\infty( \cdot| z ) $ (say). In the
latter case, given $ Z = z $, it holds that $ F_n ( \cdot| z ) \to F_\infty(
\cdot| z ) $ in the sense of the Glivenko--Cantelli theorem. Therefore, we
refer to $ F_\infty$ as the limiting e.c.d.f. [of the $p$-values
$p_i(z)$]. In view of (\ref{p_value}), we get $ F_\infty( t | z ) = P (
p_i(z) \leq t ) = 1- P(W_T(g(X_i,z)) < 1-t) = 1- P(g(X_i,z) < W_T^{-1}
(1-t) ) = 1- P(X_i < g_1( W_T^{-1} (1-t),z) ) $, hence, since $W_X$ is
assumed to be continuous,
%
%e2 ###
%
\begin{equation} \label{F_infty}
F_\infty( t | z ) = 1 - W_X \bigl( g_1 \bigl( W_T^{-1} (
1-t),z\bigr)\bigr),\qquad  t \in( 0 , 1).
\end{equation}
\end{rem}

For the sake of simplicity, it will be assumed that the model implies
that $ F_\infty( t | z ) $ is continuous in $ t \in[ 0 , 1 ] $ and
differentiable from the right at $ t = 0 $ with \mbox{$ F_\infty( 0 | z
) = 0 $} for all $ z \in\mathcal{Z}$.

In the case that a proportion $ \zeta_n =n_0 / n $ of hypotheses is
true and the rest is false, that is, $ n_0 $ hypotheses are true and $
n_1 = n - n_0$ hypotheses are false, we make the following additional
assumption in order to avoid laborious limiting considerations as $
\vartheta_i \rightarrow\infty$ under $K_i$. It will be assumed that
under an alternative $ K_i \dvtx  \vartheta_i > 0 $, the parameter
value $ \vartheta_i = \infty$ is possible. Moreover, for $ \vartheta_i
= \infty$, it will be assumed that the $p$-value $p_i$ has a Dirac
distribution with point mass in $0$. We refer to this situation as the
D-EX$(\zeta_n)$ model. As briefly mentioned in the introduction, under
suitable assumptions, EER becomes and FDR seems to become largest if $
\vartheta_i \rightarrow\infty$ for all $ i $ with $ \vartheta_i \in K_i$.
In order to calculate upper bounds for EER and FDR, we therefore
restrict attention to the D-EX$(\zeta_n)$ model which rarely (never)
appears in practical applications.
If one is interested in EER and FDR under other parameter
configurations, one may put a prior on the $\vartheta_i$'s under
alternatives $K_i$, which results in a mixture model as considered, for
example, in \cite{r12} or \cite{r22}. This makes things slightly more
complex and will not be considered in this  paper. In the
\mbox{D-EX$(\zeta_n)$} model, the e.c.d.f. of the $p$-values will be denoted by
$ F_n ( \cdot| z , \zeta _n ) $.

The following two examples fit in the D-EX$(\zeta_n)$ model and will be
studied in more detail in Sections \ref{sec4}~and~\ref{sec5},
respectively.

\begin{example} \label{Example_2.1}
Let $ X_i \sim N ( 0 , 1 ) $, $ i = 0 , \ldots, n $, be independent
standard normal random variables and let $ T_i = \vartheta_i + \sqrt{
\overline{\rho}} X_i - \sqrt{ \rho }X_0 $ with $ \vartheta_i \geq0 $, $
i = 1 , \ldots, n $, where $ \rho\in( 0 , 1 ) $ is assumed to be known
and $ \overline {\rho}= 1 - \rho$.
Then $ T = ( T_1, \ldots, T_n ) $ is multivariate normally distributed
with mean vector $ \vartheta= ( \vartheta_1 , \ldots, \vartheta_n ) $,
$ \operatorname{Var}[ T_i ] = 1 $ for $ i = 1 , \ldots, n $, and $ \operatorname{Cov} (
T_i , T_j ) = \rho$ for $ 1 \leq i \not= j \leq n $.
Consider the multiple testing problem $ H_i \dvtx  \vartheta_i = 0 $
versus $ K_i \dvtx  \vartheta_i > 0 $, $ i = 1 , \ldots, n$.
This setup includes the well-known many-one multiple comparisons
problem with underlying balanced design.
%which usually reads as follows.
%Let $ \overline{Y}_{i \cdot} \sim N ( \nu_i , \sigma^2 / m_i ) $, $ i
%= 0 , \ldots, n $,
%denote independently normally
%distributed sample means with $ \sigma^2 > 0 $ (known), $ m_1 = \cdots
%= m_n $
%and $ \nu_i \geq\nu_0 $ for $ i = 1 , \ldots, n $.
%Suppose one is interested in testing $ \tilde{H}_i : \nu_i = \nu_0 $
%versus $ \tilde{K}_i : \nu_i > %\nu_0 $ for
%$ i = 1 , \ldots, n $ by using the test statistics
%$ W_i = ( 1 / m_0 + 1 / m_1 )^{-1/2} ( \overline{Y}_{i \cdot} -
%$ i = 1 , \ldots, n $.
%Then $ \E[ W_i ] = ( 1 / m_0 + 1 / m_1 )^{-1/2} ( \nu_i - \nu_0 ) /
%$ \Var[ W_i ] = 1 $ and $ \Cov( W_i , W_j ) = \sqrt{ m / ( m + m_0 )
%} = \rho$ (say).
%
For $ \rho\in( 0 , 1 ) $, the distribution of $T$ is MTP$_2$ so that the
Benjamini--Hochberg bound applies; cf. \cite{r2}~or~\cite{r14}. Note
that $ Z $ is replaced by $ X_0 $ and $W_X=W_{X_0}=W_T=\Phi$, where $
\Phi$ denotes the c.d.f. of the standard normal distribution. Suitable
$p$-values for testing the $ H_i $'s are given by $ p_i = p_i ( x_0 ) = 1
- \Phi( \vartheta_i + \sqrt{\overline{\rho }} x_i- \sqrt{\rho} x_0 ) $,
$ i = 1 , \ldots, n $. Again, we add $ \vartheta_i = \infty$ to the
model such that $ p_i = 0$ a.s.\ if $ \vartheta_i = \infty$, $i = 1 ,
\ldots, n $.
We denote this D-EX$(\zeta_n)$ model by D-EX-N$(\zeta_n)$.
\end{example}

\begin{example} \label{Example_2.2}
Let $ X_i \sim N ( \vartheta_i , \sigma^2 ) $, $ i = 1 , \ldots, n $,
be independent normal random variables and let $ \nu S^2 / \sigma^2
\sim\chi^2_\nu$ be independent of the $X_i$'s. Without loss of
generality, we assume $ \sigma^2 = 1 $ and the c.d.f. of $ \sqrt{\nu}
S$ will be denoted by $F_{\chi_\nu}$. Again we consider the multiple
testing problem $ H_i \dvtx  \vartheta_i = 0 $ versus $ K_i \dvtx
\vartheta_i > 0 $, $ i = 1 , \ldots, n$. Let $ T_i = X_i / S $, $ i =
1, \ldots, n $.
Then $ T= ( T_1, \ldots, T_n ) $ has a multivariate equi-correlated
$t$-distribution. The c.d.f. (p.d.f.) of a univariate (central)
\mbox{$t$-distribution} will be denoted by $F_{t_\nu}$ ($f_{t_\nu}$) and a $
\beta$-quantile of the $t_\nu$-distribution will be denoted by
$t_{\nu,\beta}$. Here $ Z $ is replaced by $ S $, $W_X=\Phi$, $W_{S}(s)
=F_{\chi_\nu}(s/\sqrt{\nu})$ and \mbox{$W_T=F_{t_\nu}$}. Suitable $p$-values (as a function of
$ s $) are defined by $ p_i ( s ) = 1 - F_{t_\nu} ( x_i / s ) $. Again
we add $ \vartheta_i = \infty$ to the model such that $ p_i = 0$ a.s.
if $ \vartheta_i = \infty$.
We denote the corresponding D-EX$(\zeta_n)$ model by D-EX-t$(\zeta_n)$.
It is outlined in \cite{r2} by employing PRDS arguments that the
Benjamini--Hochberg bound applies in this model for $ \alpha \in( 0,1/2) $.
\end{example}

The following obvious lemma gives explicit expressions for $ F_\infty$
(as a consequence of the Glivenko--Cantelli theorem, cf. (\ref{F_infty}) in
Remark~\ref{rem2_1}) and characterizes crossings with the Simes-line in
the D-EX$(\zeta_n)$ model.

\begin{lem} \label{lem3}
Given D-EX$(\zeta_n)$ with $ \lim_{n \to\infty} \zeta_n = \zeta\in( 0 ,
1 ] $, the limiting e.c.d.f. of the $p$-values is given by
\[
F_\infty( t | z , \zeta) = (1-\zeta) + \zeta\bigl(1 - W_X \bigl( g_1 \bigl(
W_T^{-1} ( 1-t),z\bigr)\bigr)\bigr) , \qquad t \in( 0 , 1 ).
\]
Moreover, $ F_\infty$ crosses (or contacts) the Simes-line, that is, $
F_\infty( t | z , \zeta) = t / \alpha$ for some $ t \in(
\alpha(1-\zeta) , \alpha) $ iff
$ W_X^{-1} ( (1 - t / \alpha)/\zeta) = g_1 ( W_T^{-1} ( 1-t),z)$. If $
F_\infty( t | z ) $ is strictly decreasing in $ z $ for all $ t \in(
\alpha(1-\zeta) , \alpha) $ and if $ F_\infty( t | z , \zeta) $ $= t /
\alpha$ for some $ t^* \in( \alpha(1-\zeta) , \alpha) $, then
\[
z = z( t^* | \zeta) = g_2 \bigl( W_X^{-1} \bigl( (1 - t^* / \alpha)/\zeta\bigr),
W_T^{-1} ( 1 - t^* ) \bigr).
\]
\end{lem}

Note that $F_\infty( t | z ) = F_\infty( t | z , 1 ) $. Analogously, we
set $ z(t)=z(t|1)$.

%f1 ###
%
\begin{figure}[b] % figuur 1

\includegraphics{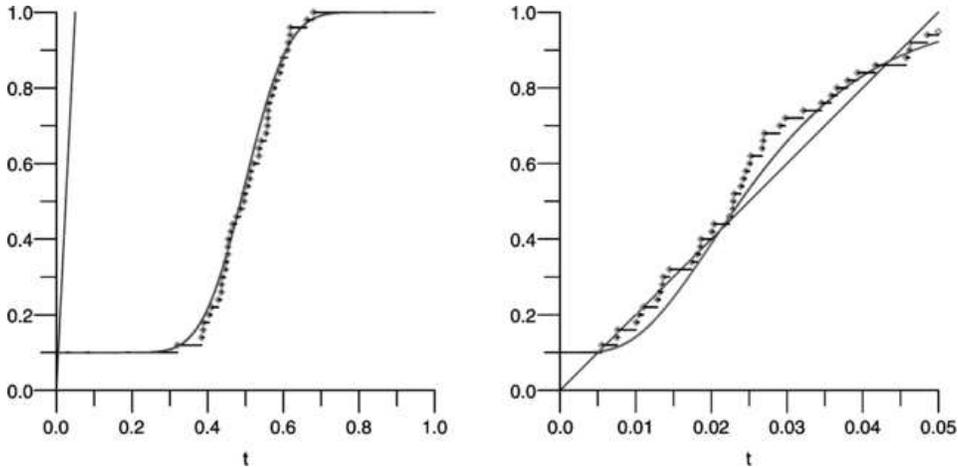}

\caption{The Simes-line for $ \alpha=0.05$ and two
realizations of the e.c.d.f. $ F_n (\cdot|x_0) $ together with $
F_\infty(\cdot|x_0) $ in the D-EX-N$(\zeta_n)$ model for $ n = 50 $, $
\zeta_n = 0.9 $, $ \rho =0.95 $ and $x_0=0.0 $ (left picture with
$t\in[0,1]$), $x_0=-2.0 $ (right picture with $t\in[0,0.05]$).}\label{figure_e.c.d.f._Simes}
\end{figure}

Figure \ref{figure_e.c.d.f._Simes} illustrates the enormous impact of
the disturbance variable and a large correlation in the
D-EX-N$(\zeta_n)$ model on the LCP determining the number of
rejections. In this example, for $ x_0 = 0.0 $ only the (totally) false
hypotheses are rejected, while for $ x_0 = -2.0 $ we obtain $38$ false
rejections.

\begin{rem}
Under the assumptions of Lemma \ref{lem3}, the Glivenko--Cantelli theorem yields
\[
\lim_{ n \to\infty} \sup_{t \in[ 0 , 1 ] } \big| F_n ( t | z , \zeta _n ) -
F_\infty( t | z , \zeta) \big| = 0 \qquad\mbox{almost surely for all } z
\in\mathcal{Z}.
\]
Moreover,
\[
\mathbb{E} [ F_{\infty} ( t | Z , \zeta) ] = \int F_{\infty} ( t | z ,
\zeta) \,dP^Z (z) = 1- \zeta+ \zeta t \qquad\mbox{for all } t \in[ 0 , 1 ].
\]
\end{rem}

For $ \lim_{n \to\infty} \zeta_n = \zeta\in( 0 , 1 ] $, define
%e3 ###
%
\begin{equation} \label{LCP}
t ( z | \zeta) = \sup\{ t \in[ \alpha( 1 - \zeta) , \alpha] \dvtx
F_\infty( t | z , \zeta) = t / \alpha\}.
\end{equation}

If there exists an $ \epsilon>0 $ such that $ F_\infty( t | z , \zeta)
> t/\alpha$ for all $ t \in[ t( z | \zeta)-\epsilon,t( z | \zeta)) $
and $ F_\infty( t | z , \zeta) < t / \alpha$ for all $ t \in( t( z |
\zeta) ,t( z | \zeta)+\epsilon] $, then $ t( z | \zeta) $ will be
called the largest crossing point (LCP) of $ F_\infty( \cdot| z ,
\zeta) $ and the Simes-line. The set of LCP's will be denoted by $
C_\zeta$. Moreover, set $ D_\zeta= \{ z \in\mathcal{Z}\dvtx  t( z |
\zeta) \in C_\zeta\} $. Note that there may be some tangent points
(TP's) $t( z | \zeta)$ defined by (\ref{LCP}) with $ F_\infty( t | z ,
\zeta) \leq t / \alpha$ in a neighborhood of $t( z | \zeta)$. However,
it will be assumed that $P^Z(D_\zeta)=1$. In practical examples, $
C_\zeta$ is a finite union of intervals. For $ \zeta\in( 0 , 1 ) $, we
always have a well defined LCP or TP $t( z | \zeta) \geq\alpha( 1 -
\zeta) > 0 $. For $ \zeta=1 $, the LCP may be $0$ for a large set of
$z$-values, which makes the calculation of the limiting EER and FDR
subtler.

In the following we make use of the notation
\begin{eqnarray*}
\operatorname{FDR}_n ( \zeta_n | z ) & = & \mathbb{E}\biggl[ \frac{V_n}{
R_n \vee1 } \Big| Z = z \biggr] ,\qquad
\operatorname{FDR}_n ( \zeta_n ) = \mathbb{E}\biggl[ \frac{V_n}{ R_n \vee
1 } \biggr],
\\
\operatorname{FDR}_\infty( \zeta| z ) & = & \lim_{n \to\infty}
\operatorname{FDR}_n ( \zeta_n | z ) ,\qquad\hspace{4pt}
\operatorname{FDR}_\infty( \zeta) = \lim_{n \to\infty}
\operatorname{FDR}_n ( \zeta_n ),
\end{eqnarray*}
and the corresponding expressions for EER. Moreover, the notation
$V_n(z)$, $R_n(z)$ will be used if $Z=z$ is given.

%s2.1 ###
\subsection{All LCP\textup{'}s greater than zero}
We first consider the case $ \zeta\in( 0 , 1 ) $.
\begin{thm} \label{thm3_1}
Given D-EX$(\zeta_n)$ with $ \lim_{n \to\infty} \zeta_n = \zeta\in(0,1)
$,  for all~$ z \in D_\zeta$
%
%e5 ###
%e4 ###
%
\begin{eqnarray}
\lim_{ n \to\infty} \frac{V_n(z)}{n} & = & \frac{t( z | \zeta
)}{\alpha} - ( 1 - \zeta)\qquad \mbox{a.s.},
\label{vn2} \\
\lim_{ n \to\infty} \frac{V_n(z)}{R_n(z) \vee1} & = &
1- \frac{ \alpha( 1 - \zeta)}{t( z | \zeta)}\qquad\hspace{15pt} \mbox{a.s.}
\label{vnrn2}
\end{eqnarray}
\end{thm}

\begin{pf}
With a similar technique as in the proof of Lemma A.2 in \cite{r9}, it
can be shown that the proportion of rejected hypotheses $ R_n ( z ) / n
$ converges almost surely to $ t( z | \zeta)/\alpha$. This fact
immediately implies (\ref {vn2}) and (\ref{vnrn2}).
\end{pf}

\begin{rem}
Under the assumptions of Theorem \ref{thm3_1},
%e7 ###
%e6 ###
%
\begin{eqnarray}
\operatorname{EER}_\infty( \zeta| z )
& = & \mathbb{E}\biggl[ \lim_{ n \to\infty} \frac{V_n(z)}{n} \biggr]
=\frac{t( z | \zeta)}{\alpha} - ( 1 - \zeta) ,
\label{vn2a} \\
\operatorname{FDR}_\infty( \zeta| z )
& = &\mathbb{E}\biggl[ \lim_{ n \to\infty} \frac{V_n(z)}{R_n(z) \vee1}
\biggr]=1- \frac{ \alpha( 1 - \zeta)}{t( z | \zeta)}.
\label{vnrn2a}
\end{eqnarray}
In view of the general assumption $P^Z ( D_\zeta) = 1 $, $z$ can be
replaced by $Z$ in (\ref{vn2}) and (\ref{vnrn2}).
\end{rem}

It remains to calculate
$ \operatorname{EER}_\infty( \zeta) $ and
$\operatorname{FDR}_\infty( \zeta) $.
This may be done in two ways. The first is to integrate (\ref{vn2}) and
(\ref{vnrn2}) with respect to $ P^Z $. In this case the main problem is
the computation of $ t ( z | \zeta) $, which can be cumbersome. In
general, $ t ( z | \zeta) $ cannot be determined explicitly.
The second possibility seems more convenient and is summarized in the
following theorem.

\begin{thm} \label{thm3_2}
Under the assumptions of Theorem \textup{\ref{thm3_1}}, suppose that $
F_\infty( t | z ) $ is strictly decreasing in $ z $ for $ t \in (
\alpha( 1 - \zeta) , \alpha) $. Let $ C_{\zeta,1} = \{
t/\alpha-1+\zeta\dvtx  t \in{C}_\zeta\} $, $ C_{\zeta,2} = \{ 1-
\alpha( 1 - \zeta) / t \dvtx  t \in{C}_\zeta\} $ and denote the c.d.f.
of\break $ \lim_{ n \to\infty} {V_n(Z)} / {n} $ and $\lim_{ n \to\infty}
{V_n(Z)}/ ({R_n(Z) \vee1)} $ by $G_{\zeta,1} $ and $G_{\zeta,2}$,
respectively. Then
%e9 ###
%e8 ###
%
\begin{eqnarray}
G_{\zeta,1} ( u ) &=& 1 - W_Z \bigl( z \bigl( \alpha( u + 1 -
\zeta) | \zeta\bigr) \bigr)\qquad \mbox{for }
u \in C_{\zeta,1} \cap( 0, \zeta),
\label{rep1} \\
G_{\zeta,2} ( u ) &=& 1 - W_Z \biggl( z \biggl( \frac{\alpha( 1 -
\zeta)}{1-u} \Big| \zeta\biggr) \biggr)\qquad\mbox{for } u \in
C_{\zeta,2} \cap( 0, \zeta) , \label{rep2}
\end{eqnarray}
hence, $\operatorname{EER}_\infty$ and $\operatorname{FDR}_\infty$ can
be computed via
%& = & \int_{C_{\zeta,1}} u d G_{\zeta,1} ( u ) \label{EER1} , \\
%
\[
\operatorname{EER}_\infty( \zeta) = \int_{C_{\zeta,1}} u\,dG_{\zeta,1}
(u) \quad\mbox{and}\quad
\operatorname{FDR}_\infty( \zeta) = \int_{C_{\zeta,2}} u \,dG_{\zeta,2} ( u ).
\]
\end{thm}

\begin{pf}
Let $ \zeta\in( 0 , 1 ) $ and $ t \in{C}_\zeta\cap( \alpha( 1 - \zeta),
\alpha) $. From (\ref{vn2}) in Theorem \ref{thm3_1}, we get
\[
\biggl\{
z \in D_\zeta\dvtx  \lim_{n \to\infty} \frac{V_n(z)}{n} > \frac
{t}{\alpha} - (1-\zeta) \mbox{ a.s. }
\biggr\}
=\{ z \in D_\zeta\dvtx  z < z( t | \zeta) \}.
\]
Hence, the substitution $ u = t / \alpha- (1-\zeta) $ yields that for
all $ u \in C_{\zeta,1} \cap( 0, \zeta)$
\begin{eqnarray*}
W_Z \bigl( z \bigl( \alpha( u + 1 - \zeta) | \zeta\bigr) \bigr) &=& P^Z \biggl( \biggl\{ z \in
D_\zeta\dvtx  \lim_{n \to\infty} \frac{V_n(z)}{n} >
\frac{t}{\alpha} - (1-\zeta) \mbox{ a.s. } \biggr\} \biggr)
\\
&=& 1 - G_{\zeta,1} ( u ),
\end{eqnarray*}
%
%for all $ u \in C_{\zeta,1} $.
%
which is (\ref{rep1}).
Similarly, we obtain from (\ref{vnrn2}) in Theorem \ref{thm3_1} that
\[
\lim_{ n \to\infty} \frac{V_n(z)}{R_n ( z ) \vee1} > 1-
\frac{\alpha(1-\zeta)}{t} \qquad\mbox{a.s.}\qquad
\mbox{iff }z < z( t | \zeta) .
\]
Therefore, similar arguments as before yield that
$ G_{\zeta,2} $ is given by (\ref{rep2}).
\end{pf}

The latter result is a key step for the computation of
$\operatorname{EER}_\infty( \zeta)$ and \linebreak[4]
$\operatorname{FDR}_\infty( \zeta) $. In practical examples it remains
to determine the sets $C_{\zeta,1}$ and $C_{\zeta,2}$ and to evaluate
the corresponding integrals.

%
%s2.2 ###
\subsection{Some LCP\textup{'}s equal to zero}

If an LCP is equal to zero, the behavior of the FDR heavily depends on
the gradient in zero of the c.d.f. of the $p$-value distribution. The
next two lemmas cover the finite case.

\begin{lem} \label{FDR_lem1}
Let $\alpha\in(0,1)$, $ 0 \leq\gamma\leq1/\alpha$, $n_0, n \in
\mathbb{N}$, $n_0 \leq n $ and let $\xi_1, \ldots, \xi_{n_0}$ be i.i.d.
random variables with values in $[0, 1]$ with c.d.f. $F_\xi$ satisfying
$F_\xi(t) = \gamma t$ for all $t \in[0, \alpha]$. Furthermore, let
$\xi_{n_0+1}, \ldots, \xi_{n}$ be random variables with values in $[0,
1]$, independent of $ ( \xi_j \dvtx  1 \leq j \leq n_0) $. For $c_i = i
\alpha/n $, $i=1, \ldots, n$, define $R'_n = \max\{k \leq n\dvtx
\xi_{k:n} \leq c_k \}$ and $V'_n = | \{ i \in\{ 1 , \ldots,n_0 \} \dvtx
\xi_i \leq c_{R_n'} \}|$ (with $ c_{R_n'} = - \infty$ for $ R_n' = -
\infty$). Then
%e10 ###
%
\begin{equation} \label{0-frage}
\mathbb{E}\biggl(\frac{V'_n}{R'_n \vee1 } \biggr) = \frac{n_0}{n}
\gamma\alpha.
\end{equation}
\end{lem}

\begin{pf}
For $ 1 \leq i \leq n_0$, denote the $(n-1)$-dimensional random vector
$(\xi_1, \ldots, \xi_{i-1}, \xi_{i+1}, \ldots, \xi_n)$ by $\xi^{(i)}$,
define for $1 \leq k < n$ the sets $D_k^{(i)}(\alpha) =
\{\xi^{(i)}_{k:n-1} > c_{k+1}, \ldots, \xi^{(i)}_{n-1:n-1} > c_{n}\}$
and set $D_0^{(i)}(\alpha) = \varnothing$, $D_n^{(i)}(\alpha) =
\Omega$.
%Notice that the sets $D_k^{(i)}(\alpha)$ are increasing in $k$.
Then the left-hand side of (\ref{0-frage}) (cf., e.g., Lemma 3.2 and
formula (4.4) in \cite{r14}) is equal to
\[
\frac{1}{n} \sum_{i=1}^{n_0} \mathrm{P} (\xi_i \leq c_n)
+ \sum_{i=1}^{n_0} \sum_{j=2}^{n} \biggl[\frac{\mathrm{P}(\xi_i
\leq c_{j-1})}{j-1} - \frac{\mathrm{P}(\xi_i \leq c_{j})}{j}\biggr]
\mathrm{P}\bigl(D_{j-1}^{(i)}(\alpha)\bigr).
\]
Noting that $ \mathrm{P} (\xi_i \leq c_n) = \gamma\alpha$ for all $ 1
\leq i \leq n_0 $ and $ \mathrm{P}(\xi_i \leq c_j)/j = \gamma\alpha/ n $
for all $ 1 \leq j \leq n $, the assertion follows immediately.
\end{pf}

As an application of Lemma \ref{FDR_lem1}, we insert a surprising
example.

\begin{example} \label{examp_expo}
Under the general framework of this section, suppose the $ X_i$ follow
an exponential distribution with scale parameter $ \lambda=1$ and
location parameter $ \vartheta_i $ and $ Z $ follows an exponential
distribution with scale parameter $ \lambda=1$ and location parameter
$0$, and consider the model $ T_i = g(X_i , Z) = X_i - Z$, $ i = 1 ,
\ldots,n $. Under $ H_i\dvtx  \vartheta_i = 0 $, the c.d.f. of $T_i$ is
given by
%W_T(t) =
% \end{cases}
$W_T(t) = \exp(t)/2 $ for $ t \leq0 $ and $W_T(t) = 1 -\exp(-t)/2$ for
$ t > 0 $, while the $p$-values (as functions of the observed $z$-value)
are given by $ p_i(z) = 1-W_T ( x_i - z ) $, $ i = 1 , \ldots,n $. This
results in
\[
F_{\infty}(t|z) =
\cases{
2 \exp(-z) t, &\quad\mbox{if }$0 \leq t \leq1/2$,\cr
\exp(-z) (2-2t)^{-1},&\quad\mbox{if }$1/2 < t \leq u(z)$,\cr
1, &\quad\mbox{if }$u(z) < t \leq1$,}
\]
with $u(z) = 1 - \exp(-z)/2$. For convenience, we restrict attention to
the case \mbox{$ \alpha\leq1/2 $}. In order to apply Lemma \ref{FDR_lem1}, set
$ F_\xi( t )= F_{\infty}(t|z)$ and note that $ p_i ( z) $ has c.d.f. $
F_\xi$ if $H_i$ is true. Therefore, supposing that $n_0$ hypotheses are
true and $n_1=n-n_0$ are false with fixed but arbitrary $ \vartheta_i >
0 $, we obtain with $ \gamma( z ) = 2 \exp(-z)$ and $ \zeta_n = n_0 / n
$ that $ \operatorname{FDR}_n ( \zeta_n | z ) = \zeta_n \alpha\gamma( z
) \mbox{ for all } z > 0 . $ Integrating with respect to $ P^Z $
finally results in
\[
\operatorname{FDR}_n ( \zeta_n ) = \zeta_n \alpha\int\gamma(z) \,dP^Z
(z) = \zeta_n \alpha.
\]
It may be astonishing that the Benjamini--Hochberg upper bound for the
FDR is attained for all parameter configurations although the $ T_i $'s
are dependent. Notice that the $\mbox{MTP}_2$ property holds in this
setting so that the Benjamini--Hochberg bound for the FDR applies. This
is a consequence of Propositions 3.7 and 3.8 in \cite{karlnott-1980},
because the p.d.f. of the $\operatorname{Exp}(\lambda)$ distribution is
$\mbox{PF}_2$ for any $\lambda> 0$; see \cite{karlin-1968}, page~30.
\end{example}

The next result extends Lemma \ref{FDR_lem1} and is a helpful tool in
the case that LCP's are in $0$.
\begin{lem} \label{lem3_4}
Under the assumptions of Lemma \textup{\ref{FDR_lem1}} but only
supposing that $F_\xi(t) = \gamma t $ for all $t \in[0, t^*]$ for some
$ t^* \in(0, \alpha]$, let $ A_n(t^*) =\{F_n(t) < t/\alpha\,\forall t
\in(t^*, \alpha]\} $, where $F_n$ denotes the e.c.d.f. of $\xi_1,
\ldots,\xi_n$. Then, setting $r = \max\{i \in\mathbb{N}_0\dvtx  i
\alpha/ n \leq t^* \}$,
%
%e11 ###
%
\begin{equation} \label{restringiert}
\mathbb{E}\biggl(\frac{V'_n}{R'_n \vee1} \mathbf{1}_{A_n(t^*)} \biggr) =
\frac{n_0}{n} \gamma\alpha\mathrm{P}\bigl(D_r^{(1)}(\alpha)\bigr).
\end{equation}
\end{lem}

\begin{pf}
It is clear that $A_n(t^*) = \{R'_n \leq r\}$; hence, for $r>0$, the
left-hand side of (\ref{restringiert}) is now equal to
\[
\frac{1}{r} \sum_{i=1}^{n_0} \mathrm{P}(\xi_i \leq c_r)
\mathrm{P}\bigl(D_r^{(i)}(\alpha)\bigr) + \sum_{i=1}^{n_0} \sum_{j=2}^{r}
\biggl[\frac{\mathrm{P} (\xi_i \leq c_{j-1})}{j-1} - \frac{\mathrm{P}(\xi_i \leq
c_{j})}{j}\biggr] \mathrm{P}\bigl(D_{j-1}^{(i)}(\alpha)\bigr).
\]
The assertion follows  similarly as in the proof of Lemma \ref{FDR_lem1}.
\end{pf}
The following theorem, the proof of which is in the
\hyperref[app]{Appendix}, is an important step for the understanding of
the limiting behavior of both EER (or ENE) and FDR given a fixed value
$Z=z$ such that the LCP is in $0$.

\begin{thm} \label{thm2_3}
Given D-EX$(\zeta_n)$ with $ \lim_{n \to\infty} \zeta_n = 1 $, let $ z
\in\mathcal{Z}$ be such that $ F_\infty(t|z) < t/\alpha $ for all $ t
\in(0,\alpha] $. Setting $ \gamma(z) = \lim_{t \to0^+} F_\infty(t|z) /
t , $ it holds that
%
%e13 ###
%e12 ###
%
\begin{eqnarray}
\operatorname{EER}_\infty(1|z) &=& 0 ,
\label{vn3} \\
\operatorname{FDR}_\infty( 1 |z) & = & \alpha\gamma( z ).
\label{vnrn4}
\end{eqnarray}
\end{thm}

\begin{rem}
In \cite{r10} the distribution and expectation of $V_n$ were computed
for uniform $p$-values under the assumption that all hypotheses are true.
Assuming $ \zeta_n = 1 $ for all $ n \in\mathbb{N}$, the nesting method
in the proof of (\ref{vnrn4}) together with the technique in \cite{r10}
can be used to prove
\[
\lim_{ n \to\infty} \mathbb{E}[ V_n(z) ] =
\cases{
\displaystyle\frac{\alpha\gamma( z) }{( 1 - \alpha\gamma( z))^2},&\quad $\gamma( z ) < 1/\alpha$,
\cr
\infty,&\quad $\gamma( z ) = 1/\alpha$.}
\]
It is important to note that this formula is only valid for $ \zeta_n =
1 $. If $ n_1 $ tends to infinity with $\lim_{n \to\infty} n_1 / n = 0
$ and $ \gamma( z ) > 0 $, then $ \lim_{ n \to\infty} \mathbb{E}[
V_n(z) ] = \infty$.
\end{rem}

To complete the picture for $ \zeta=1 $, the next theorem puts things
together.

\begin{thm} \label{thm2_4}
Given D-EX$(\zeta_n)$ with $ \lim_{n \to\infty} \zeta_n = 1 $, suppose
that $ F_\infty( t | z ) $ is strictly decreasing in $ z $ for $ t \in(
0 , \alpha) $.
Moreover, let $G_{1,1}$ and $ C_{1,1} $ be defined as in
Theorem~\textup{\ref{thm3_2}} and let $ E_0 = \{ z \in\mathcal{Z}\dvtx
t(z|1) = 0 \} $ and $ E_1=\mathcal{Z}\setminus E_0 $. Then
%
%e15 ###
%e14 ###
%
\begin{eqnarray}
\operatorname{EER}_\infty( 1 )
& = & \int_{C_{1,1}} u \,d G_{1,1} ( u ) ,
\label{EER10} \\
\operatorname{FDR}_\infty( 1 ) & = & P^Z(E_1) + \alpha\int_{E_0}
\gamma(z) \,d P^Z (z) . \label{FDR10}
\end{eqnarray}
\end{thm}

\begin{pf}
Using the disjoint decomposition $\mathcal{Z}= E_0 + E_1$, we obtain
\begin{eqnarray*}
\operatorname{EER}_\infty( 1 ) &= &\lim_{n \to\infty} \int_{\mathcal{Z}} \frac{V_n(z)}{n} \,dP^Z(z)
\\
&= &\int_{E_0} \lim_{n \to\infty} \frac{V_n(z)}{n} \,dP^Z(z) +
\int_{E_1} \lim_{n \to\infty} \frac{V_n(z)}{n} \,P^Z(z)
\\
&= &A_1 + A_2 \mbox{ (say)}.
\end{eqnarray*}
Now, Theorem \ref{thm2_3} immediately yields $A_1 = 0$, and in analogy
to the arguments in the proof of Theorem \ref{thm3_2}, we get that
$A_2= \int_{C_{1,1} \setminus\{ 0 \} } u\,d G_{1,1} ( u )$. Therefore,
(\ref{EER10}) is proven. Applying the same decomposition (together with
Theorem \ref{thm2_3}) to $\operatorname{FDR}_\infty( 1 )$ and observing
that $\lim_{n \to \infty} V_n(z) / (R_n(z) \vee1) = 1$ if $z \in E_1$
[similar to (\ref {vnrn2}) with $\zeta= 1$] finally
proves~(\ref{FDR10}).
\end{pf}

%s3 ###
\section[Exchangeable normal variables (Example 2.1 continued)]{Exchangeable normal variables (Example \protect\ref{Example_2.1} continued)} \label{sec4}

In the D-EX-N$(\zeta_n)$ model, assuming that the proportion $ \zeta_n
$ of true null hypotheses tends to 1, we obtain from Lemma \ref{lem3}
that the limiting e.c.d.f. of the $ p_i $'s given $ X_0 = x_0 $ is
given by
\[
F_{\infty} ( t | x_0 ) = 1 - \Phi\bigl( \Phi^{-1}(1-t) /
\sqrt{\overline{\rho}} + \sqrt{ \rho / \overline{\rho}} x_0 \bigr)
\qquad\mbox{for all } t \in( 0 , 1 ),
\]
and $ F_{\infty} ( 0 | x_0 ) = 1 - F_{\infty} ( 1 | x_0 ) = 0 $. Note
that $ F_\infty( t | x_0 ) = P ( \sqrt{\overline{\rho}} X -\break \sqrt{\rho}
x_0 > u_{1-t} )$, where $ X $ denotes a standard normal variate and $
u_{\alpha} $ denotes\break the corresponding $ \alpha$-quantile. Moreover, it
is $ \lim_{t \downarrow0} (\partial/ \partial t )F_{\infty} ( t | x_0 )
=\break \lim_{t \uparrow1} (\partial/ \partial t ) F_{\infty} ( t | x_0 ) = 0
$, and $ F_{\infty} ( \cdot| x_0 ) $ is convex for $ 0 \leq t \leq\Phi
( x_0 / \sqrt{ \rho}) $ and concave for $ \Phi( x_0 / \sqrt{ \rho})
\leq t \leq1$. Furthermore, $ F_{\infty} ( t | x_0 ) $ is strictly
decreasing in $ x_0 $ for $ t \in( 0 , 1 ) $ and
$\lim_{\rho\downarrow0} F_{\infty} ( t | x_0 ) = t $.

Assuming that $ \lim_{ n \to\infty} \zeta_n = \zeta\in( 0 , 1 ] $, the
limiting e.c.d.f. is given by
\[
F_{\infty} ( t | x_0 , \zeta) = (1-\zeta)
+ \zeta F_{\infty} ( t | x_0 ).
\]
Hence, for $ \zeta\in(0,1]$ and given $t \in(0, \alpha)$,
$
F_{\infty} ( t | x_0 , \zeta) = t / \alpha
$
iff
%e16 ###
%
\begin{equation} \label{x_0}
x_0 = x_0 ( t | \zeta) = \sqrt{ \overline{\rho}/ \rho} \Phi
^{-1}\bigl((1-t/\alpha)/\zeta\bigr) - \Phi^{-1}(1-t) /\sqrt{\rho}.
\end{equation}

For $ \zeta\in(0,1) $, we get $ \lim_{t \downarrow\alpha( 1 - \zeta ) }
x_0 ( t | \zeta) = + \infty $ and $ \lim_{t \uparrow\alpha} x_0 ( t |
\zeta) = - \infty $. Moreover, $F_\infty( \cdot| x_0, \zeta) $ starts
above the Simes-line so that there is at least one CP in $ (0,1)$. In
fact, there may be one, two or three points of intersection in $
(0,1)$. For $ \zeta=1 $, we get in contrast to $ \zeta\in(0,1) $ that $
\lim_{t \downarrow0 } x_0 ( t | \zeta) = \lim_{t \uparrow\alpha } x_0 (
t | \zeta) = - \infty $.
The limiting e.c.d.f. $ F_\infty( \cdot| x_0)=F_\infty( \cdot| x_0, 1 )
$ starting with $ F_\infty( 0 | x_0 ) =0 $ may have no, one or two CP's
in $ (0,1)$.

In order to determine the set of LCP's, the following derivations are
helpful. Let $ u = \Phi^{-1} ( 1 - t ) $ and let
%e17 ###
%
\begin{equation} \label{distanz}\hspace*{20pt}
d ( u | x_0 , \zeta) = (1-\zeta) + \zeta\bigl(1 - \Phi\bigl(u /
\sqrt{\overline{\rho}} + \sqrt{ \rho/ \overline{\rho}} x_0
\bigr)\bigr) - \bigl( 1 - \Phi( u ) \bigr) / \alpha
\end{equation}
denote the distance between the transformed $F_\infty$-curve and the
transformed Simes-line. Then the conditions
%e19 ###
%e18 ###
%
\begin{eqnarray}
d ( u | x_0 , \zeta) & = & 0,
\label{BP_cond_1} \\
\frac{\partial}{\partial u} d ( u | x_0 , \zeta) & = & 0
\label{BP_cond_2}
\end{eqnarray}
are necessary and sufficient for a TP ($ F_\infty$ touches the Simes-line).
Note that condition~(\ref{BP_cond_2}) is equivalent to
%
%e20 ###
%
\begin{equation} \label{BP_cond_3}
u \in\bigl\{ u_{1,2}( x_0 ) = - x_0 / \sqrt{ \rho} \pm\sqrt{
\overline{\rho}/ \rho}
\sqrt{ x_0^2 - 2 \ln ( \sqrt{ \overline{\rho}}/ ( \alpha\zeta) )} \bigr\}.
\end{equation}
If there exists a real solution $ u^*$ of (\ref{BP_cond_1}) and
$u^* = u_2 (x_0) = - x_0 / \sqrt{ \rho} - \sqrt{ \overline{\rho}/
\rho}\times\break
\sqrt{ x_0^2 - 2 \ln ( \sqrt{ \overline{\rho}}/ ( \alpha\zeta) ) } $
for given values of $x_{0}, \rho, \zeta$, then we define $ t_2 =\break 1 -
\Phi(u^*)$.
If such a solution $u^*$ exists in case of $ \zeta\in(0,1) $,
define $t_1$ as the smaller solution of $ F_\infty( t | x_0 , \zeta)
= t / \alpha$. Then the set of LCP's is given by $ C_\zeta= ( \alpha(
1-\zeta) , t_1 ) \cup( t_2 , \alpha) $.
Note that for $ \zeta= 1 $ there exists a unique TP such that the set
of LCP's is given by $ C_\zeta= \{ 0 \} \cup( t_2 , \alpha) $.
Furthermore, for $ \zeta\in(0,1) $, there may be no such TP. In the
latter case, formally interpreted as $t_1=t_2$, we have $ C_\zeta= (
\alpha( 1-\zeta) , \alpha) $.
%For example, if $ \rhoq\geq(\alpha\zeta)^2 $,
%we have a unique CP in $ ( 0 , \alpha) $
%iff $ d ( u_2(\underline{x}_0) | \underline{x}_0 , \zeta) > 0 $ %\geq
%0 $
%for $ \underline{x}_0 = - \sqrt{ 2 \ln( \sqrt{ \rhoq}/ ( \alpha
For example, such a situation occurs in the case  $ \overline{\rho}\geq
(\alpha\zeta)^2 $ and $ \alpha\in( 0 , 1/2] $ iff $ d (
u_2(\underline{x}_0) | \underline{x}_0 , \zeta) \geq0 $ for $
\underline{x}_0 = - \sqrt{ 2 \ln ( \sqrt{ \overline{\rho}}/ (
\alpha\zeta) ) } $.

The (discontinuous) case $ t_1 < t_2 $ looks somewhat paradoxical. In
this case, depending on the observed $ x_0 $, either a small proportion
$ \pi_1 \in(( 1 - \zeta) , t_1 / \alpha) $ or a larger proportion $
\pi_2 \in( t_2 / \alpha, 1 ) $ of hypotheses will be rejected although
the distance between the corresponding $x_0$ values may be small. This
occurs, for example, for $ \alpha=0.1$, $\zeta=0.9999$.

The following two theorems give formulas for $
\operatorname{EER}_\infty$ and $\operatorname{FDR}_\infty$. The first
theorem covers $ \zeta\in(0,1) $, the second one $ \zeta= 1$. The proof
of Theorem \ref{thm4_1} can be found in the \hyperref[app]{Appendix},
while the proof of Theorem \ref{thm4_2} is a straightforward
application of Theorem \ref{thm2_4}.
\begin{thm} \label{thm4_1}
Given model D-EX-N$(\zeta_n)$ with $ \lim_{n \to\infty} \zeta_n =
\zeta\in(0,1) $, the set of LCP\textup{'}s is $ C_\zeta= ( \alpha( 1-\zeta) ,
t_1 ) \cup( t_2 , \alpha) $ for $ t_1 < t_2 $ and $ C_\zeta= ( \alpha(
1-\zeta) , \alpha) $ for $ t_1 = t_2 $ (i.e., no TP) and
\begin{eqnarray*}
\operatorname{EER}_\infty( \zeta)
& = & \frac{t_2 - t_1}{ \alpha} \Phi( x_0 ( t_1 | \zeta) )
\\
&&{} +\int_{1-\zeta}^{t_1/\alpha} \Phi( x_0 ( \alpha t | \zeta) ) \,dt
+ \int_{t_2 / \alpha}^1 \Phi( x_0 ( \alpha t | \zeta) )\, dt,
\\
\operatorname{FDR}_{\infty} ( \zeta)
& = & ( z_2 - z_1 ) \Phi\biggl(x_0 \biggl( \frac{\alpha( 1 - \zeta)}{1 - z_1 } \Big|
\zeta\biggr) \biggr)
\\
&&{} + \int_{ 0 }^{ z_1 } \Phi\biggl( x_0 \biggl( \frac{\alpha( 1 - \zeta)}{1 - z }
\Big| \zeta\biggr) \biggr)\, d z
+ \int_{ z_2 }^{ \zeta} \Phi\biggl( x_0 \biggl( \frac{\alpha( 1 - \zeta)}{1 -
z } \Big| \zeta\biggr) \biggr)\, d z,
\end{eqnarray*}
where $ z_i = 1- \alpha( 1 - \zeta) / t_ i $ , $ i = 1 , 2 $.
\end{thm}

\begin{thm} \label{thm4_2}
Given model D-EX-N$(\zeta_n)$ with $ \lim_{n \to\infty} \zeta_n =
\zeta=1 $, the set of LCP's is $ C_\zeta= \{ 0 \} \cup( t_2 , \alpha) $
and
\begin{eqnarray*}
\operatorname{EER}_\infty( 1 )
& = & t_2 \Phi( x_0 ( t_2 | \zeta) ) / \alpha
+ \int_{ t_2 / \alpha}^1\Phi( x_0 ( \alpha t | \zeta) ) \,d t ,
\\
\operatorname{FDR}_{\infty} ( 1 )
& = & \Phi( x_0 (t_2 | \zeta) ).
\end{eqnarray*}
\end{thm}

\begin{rem}
For $ \zeta= 1 $, we obtain an upper bound for $ x_0 (t_2 | \zeta) $ and
$ \operatorname{FDR}_\infty( 1 ) $, respectively, if $ \rho\leq
1-\alpha^2 $.
From the derivations before Theorem \ref{thm4_1}, we get $ x_0 (t_2 |
\zeta) \leq\underline{x}_0 = - \sqrt{ 2 \ln ( \sqrt{ \overline{\rho}}/
\alpha ) } $ and consequently, $ \operatorname{FDR}_\infty( 1 )
\leq\Phi(\underline{x}_0)$.
This is helpful for the numerical determination of $ x_0(t_2 | \zeta) $.
\end{rem}

The following interesting and maybe unexpected result, which will be
discussed in Section \ref{sec6}, concerns a discontinuity for $ \zeta=
1 $ and $ \rho\to0 $. The proof is given in the \hyperref[app]{Appendix}.

\begin{thm} \label{thm4_3}
Given model D-EX-N$(\zeta_n)$ with $ \lim_{n \to\infty} \zeta_n =
\zeta=1 $ and $ \alpha\in( 0 , 1/2 ] $,
%e21 ###
%
\begin{equation}
\lim_{\rho\to0^+} \operatorname{FDR}_{\infty} ( 1 ) = \Phi
\bigl(-\sqrt{-2 \ln(\alpha)} \bigr). \label{surprise1}
\end{equation}
\end{thm}
\begin{figure}[t]

\includegraphics{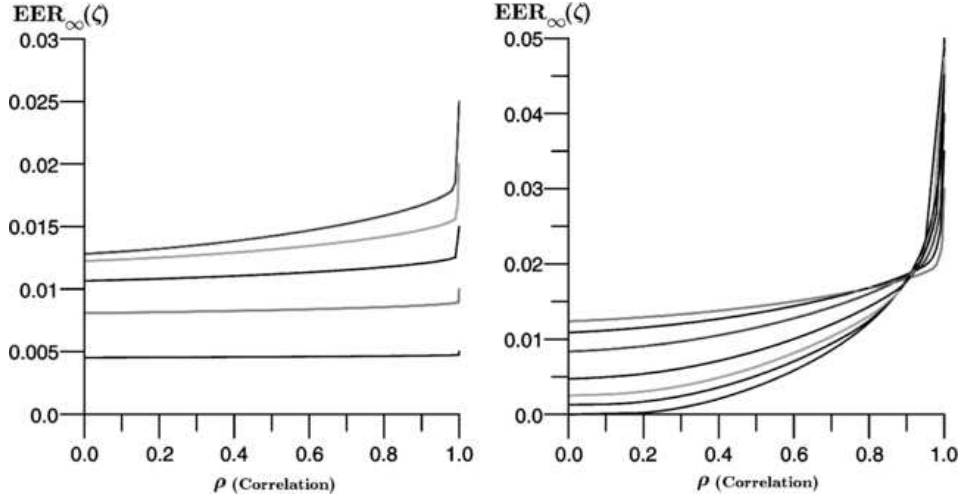}

\caption{$ \operatorname{EER}_\infty(\zeta)$ in the
D-EX-N$(\zeta_n)$ model for $\alpha=0.05$ and $\zeta=
0.1,~0.2,~0.3,~0.4,~0.5$ (left picture) and $
\zeta=0.6,~0.7,~0.8,~0.9,~0.95,~0.975,~1$ (right
picture).}\label{figure_rho_eer}
\end{figure}
\begin{figure}[b]%[H] % figuur 1

\includegraphics{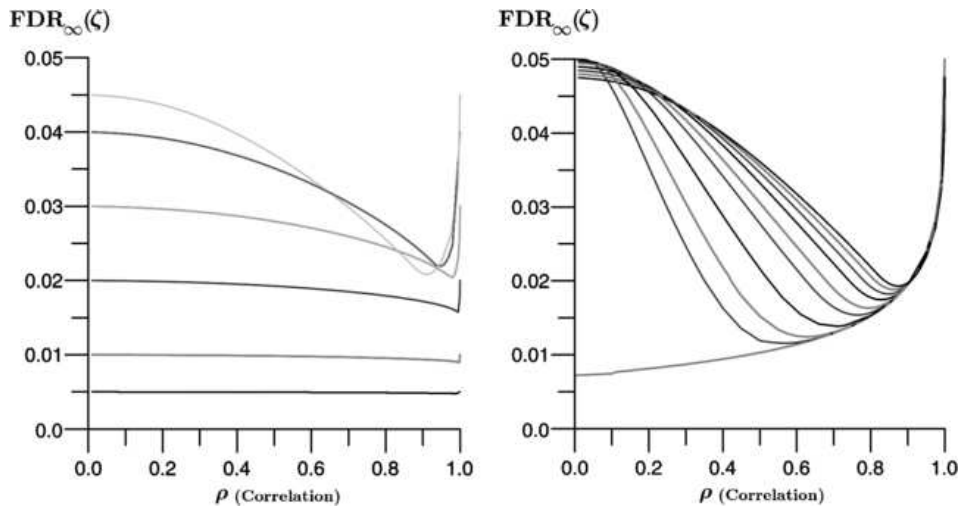}

\caption{$ \operatorname{FDR}_\infty(\zeta)$ in the
D-EX-N$(\zeta_n)$ model for $\alpha=0.05$ and $\zeta=
0.1,~0.2,~0.4,~0.6,~0.8,$ $0.9$ (left picture) and $
\zeta =0.95,~0.96,~0.97,~0.98,~0.99,~0.995,~0.999,~0.9999,~0.99999,~1$
(right picture).} \label{figure_rho}
\end{figure}

Figures \ref{figure_rho_eer} and \ref{figure_rho} display $
\operatorname{EER}_\infty( \zeta) $ and $ \operatorname {FDR}_\infty(
\zeta) $, respectively, for various values of $ \zeta$ for $
\rho\in[0,1]$. For $ \rho\to0 $, $ \operatorname{EER}_\infty( \zeta) $
tends to $ \alpha(1-\zeta)/(1-\alpha\zeta) $ as expected
(cf.~\cite{r9}) and for $ \rho\to1 $, $ \operatorname{EER}_\infty(
\zeta) $ tends to $\zeta\alpha$. Moreover, it seems that $
\operatorname{EER}_\infty( \zeta) $ is increasing in $ \rho$ with
largest values for large $\rho$ and $\zeta$. If $\rho$ is not too large
($<0.9$), $ \operatorname{EER}_\infty( \zeta) $ is largest for
$\zeta\approx1/2 $. For $ \zeta\in(0,1) $,  FDR tends to the
Benjamini--Hochberg bound for $ \rho\to0 $ and $ \rho\to1 $. For $ \rho=
1 $, we have total dependence so that $ \operatorname{FDR}_n(\zeta_n) =
\zeta_n \alpha$ in the D-EX-N$(\zeta_n)$ model. For large values of~$
\zeta$, the computation of $ \operatorname {FDR}_\infty( \zeta) $ is
extremely cumbersome. The main reason is that the TP's are very close
to 0 so that an enormous numerical accuracy is required.
Finally, it is interesting to note that for $\zeta= 1$, $
\operatorname{FDR}_\infty( 1 ) $ reflects the limiting behavior of the
true level of Simes' \cite{r17} global test for the intersection
hypothesis. Our results imply that this global test has an asymptotic
level greater than zero for all correlations $\rho\in[0, 1]$, which is a
new finding.

%s4 ###
\section[Studentized normal variables (Example 2.2 continued)]{Studentized normal variables (Example \protect\ref{Example_2.2}
continued)} \label{sec5}

In the D-EX-t$(\zeta_n)$ model with $ \lim_{n \rightarrow\infty}
\zeta_n = \zeta= 1 $, $ F_\infty( \cdot|s ) $ is given by
\[
F_\infty(t|s) = 1 - \Phi\bigl( s F_{t_\nu}^{-1} ( 1 - t ) \bigr) = 1 - \Phi(
s t_{\nu,1-t} ).
\]
Note that $F_\infty(t|s)$ is decreasing in $s$ for $ t < 1/2 $ and
increasing in $s$ for $ t > 1/2 $.
Moreover, $ (\partial/ \partial t ) F_{\infty} ( t | s ) = s \varphi( s
t_{\nu ,1-t}) / f_{t_\nu} ( t_{\nu,1-t} ) , $ hence, we get $ \lim_{t
\downarrow0} (\partial/ \partial t )F_{\infty} ( t | s ) = $ $ \lim_{t
\uparrow1} (\partial/ \partial t )F_{\infty} ( t | s ) = 0 $.
Moreover,
\[
\frac{\partial^2}{ \partial t^2 } F_{\infty} ( t | s ) > 0\qquad
\mbox{iff } - s^2 t_{\nu,1-t} < \frac{ f_{t_\nu}' ( t_{\nu,1-t} ) }{
f_{t_\nu} (t_{\nu,1-t}) }.
\]
This condition is equivalent to
\[
s^2 t_{\nu,1-t} >  \frac{ \nu+1 }{ \nu} \biggl( 1 + \frac{
t_{\nu,1-t}^2 }{\nu} \biggr)^{-1} t_{\nu,1-t}.
\]
Hence, for $ t < 1 / 2 $, $ F_\infty( t | s ) $ is convex for $ t <
\min\{ 1/2, F_{t_\nu}(- a(s,\nu) ) \} $ with $ a(s,\nu)= ( (\nu+1)/s^2
- \nu)^{1/2} $ and concave otherwise. For $ t > 1 / 2 $, $ F_\infty( t
| s ) $ is convex for $ t < \max\{ 1/2, F_{t_\nu}( a(s,\nu) ) \} $ and
concave otherwise. Notice that $ F_\infty( 1/2 | s ) = 1 / 2 $ for all
$ s>0$. As a consequence, for $ \alpha< 1/2 $, $ F_\infty$ crosses the
Simes-line at most if $ F_{t_\nu}(- a(s,\nu) ) < 1/2 $, which happens if
$ s^2 < ( \nu+ 1 ) / \nu$.

Given the D-EX-t$(\zeta_n)$ model with $ \lim_{ n \to\infty} \zeta _n =
\zeta\in( 0 , 1 ] $, the limiting e.c.d.f. is given by
\[
F_\infty(t|s, \zeta) = ( 1 - \zeta) + \zeta\bigl( 1 - \Phi\bigl( s F_{t_\nu
}^{-1} ( 1 - t ) \bigr) \bigr) .
\]
For convenience, we restrict attention to $ \alpha\in(0,1/2] $ in the
remainder of this section. For $ \zeta\in(0,1] $, we have $
F_\infty(t|s, \zeta) = t / \alpha$ iff
\[
s = s(t | \zeta)=\frac{ \Phi^{-1} ( (1 - t/\alpha)/\zeta)
}{F_{t_\nu}^{-1} ( 1-t) } ,
\]
where $ s(t | \zeta) > 0 $ iff $ t < \alpha(1-\zeta/2 )$. Therefore,
LCP's are only possible in $ [ t_u, t_o ] $ with $t_u=\alpha(1-\zeta)$
and $t_o= \alpha(1-\zeta/2)$. Notice that $ \lim_{t \downarrow t_u} s(t
| \zeta) = \lim_{t \uparrow t_o} s(t | \zeta) = 0 $ for $ \zeta= 1 $,
while $ \lim_{t \downarrow t_u} s(t | \zeta) = \infty$ and $ \lim_{t
\uparrow t_o} s(t | \zeta) = 0 $ for $ \zeta\in(0,1) $. For $
\zeta\in(0,1) $, the set $C_\zeta$ of LCP's consists of one or two
intervals denoted by $ ( \alpha( 1 - \zeta) , t_1 ) $ and $ ( t_2 ,
\alpha( 1 - \zeta/ 2 ) ) $. If there exists a TP we have $ t_1 < t_2 $
and the TP is $t_2$; otherwise $t_1=t_2$. In the case  $ \zeta= 1 $ the
existence of a TP (denoted by $t_2$) is guaranteed and the set of LCP's
is given by $C_1= \{ 0 \} \cup(t_2,\alpha/ 2 ) $. Hence, the situation
is quite similar to the D-EX-N$(\zeta_n)$ model in Section \ref{sec4}
except that there are no crossing points at all in $[t_o,\alpha]$. With
$ u = F_{t_\nu}^{-1} ( 1 - t ) $, the distance function between $
F_\infty$ and the Simes-line is defined by $ d_\nu( u | s , \zeta) =
(1-\zeta) + \zeta(1 - \Phi(s u) ) - F_{t_\nu}(-u)/\alpha. $ Necessary
and sufficient conditions for a TP ($ F_\infty$ touches the Simes-line)
are now given by
\[
d_\nu( u | s , \zeta) = 0 \quad \mbox{and}\quad
\frac{\partial}{\partial u} d_\nu( u | s , \zeta) = 0,
\]
which are equivalent to $\alpha\Phi(-s u) = F_{t_\nu}( -u)
\label{s_u_1} $ and $s \alpha\varphi( s u ) = f_{t_\nu}(u)
\label{s_u_2}$.

We summarize the behavior of EER and FDR in the following two theorems
in analogy to the results of Section \ref{sec4}.

\begin{thm} \label{thm5_1}
Given model D-EX-t$(\zeta_n)$ with $ \lim_{n \to\infty} \zeta_n =
\zeta\in(0,1) $ and $ \alpha\in(0,1/2 ] $, the set of LCP\textup{'}s is $
C_\zeta= ( \alpha( 1-\zeta) , t_1 ) \cup( t_2 , \alpha(1-\zeta/2) ) $
for $ t_1 < t_2 $ and $ C_\zeta= ( \alpha( 1-\zeta) , \alpha(1-\zeta/2)
) $ for $ t_1 = t_2 $ (i.e., no TP) and
\begin{eqnarray*}
\operatorname{EER}_\infty( \zeta)
& = & \frac{t_2 - t_1}{ \alpha} F_{\chi_\nu} \bigl( \sqrt{\nu} s ( t_1 |
\zeta) \bigr)
\\
&&{} + \int_{1-\zeta}^{t_1/\alpha} F_{\chi_\nu} \bigl( \sqrt{\nu} s ( \alpha
t | \zeta) \bigr)\, d t
\\
&&{}+ \int_{t_2 / \alpha}^{1-\zeta/2} F_{\chi_\nu} \bigl( \sqrt{\nu} s (
\alpha t | \zeta) \bigr) \,d t,
\\
\operatorname{FDR}_{\infty} (\zeta) & = & ( z_2 - z_1 ) F_{\chi_\nu} \biggl(
\sqrt{\nu} s \biggl( \frac{\alpha( 1 -
\zeta)}{1 - z_1 } \Big| \zeta\biggr) \biggr)
\\
&&{} + \int_{ 0 }^{ z_1 } F_{\chi_\nu} \biggl( \sqrt{\nu} s \biggl( \frac{\alpha
( 1 - \zeta)}{1 - z } \Big| \zeta\biggr) \biggr)\, d z
\\
&&{}+ \int_{ z_2 }^{ z_3 } F_{\chi_\nu} \biggl( \sqrt{\nu} s \biggl( \frac{\alpha
( 1 - \zeta)}{1 - z } \Big| \zeta\biggr) \biggr) \,d z,
\end{eqnarray*}
where $ z_i = 1- \alpha( 1 - \zeta) / t_ i $ , $ i = 1 , 2 $, and $
z_3=\zeta/ (2- \zeta) $.
\end{thm}

\begin{thm} \label{thm5_2}
Given model D-EX-t$(\zeta_n)$ with $ \lim_{n \to\infty} \zeta_n =
\zeta=1 $ and $ \alpha\in(0,1/2 ] $, the set of LCP\textup{'}s is $ C_\zeta= \{
0 \} \cup( t_2 , \alpha/2 ) $ and
\begin{eqnarray*}
\operatorname{EER}_\infty( 1 ) & = & t_2 F_{\chi_\nu} \bigl( \sqrt{\nu}
s ( t_2 | \zeta)\bigr) / \alpha + \int_{ t_2 / \alpha}^1 F_{\chi_\nu}
\bigl( \sqrt{\nu} s ( \alpha t | \zeta) \bigr)\,d t,
\\
\operatorname{FDR}_{\infty} ( 1 )& = &
F_{\chi_\nu} \bigl( \sqrt{\nu} s ( t_2 | \zeta) \bigr).
\end{eqnarray*}
\end{thm}

Finally, for $ \zeta= 1 $, we consider the case that the degrees of
freedom $ \nu$ tend to infinity. Heuristically, this means that the
model tends to independence.
In contrast to the normal case of the previous section, the solution is
more difficult. The reason is that one has to find suitable asymptotic
expansions for $ f_{t_\nu} $ and $ F_{t_\nu} $ given in~\cite{r4fd}.
Application of these expansions yields the following result, the proof
of which is given in the \hyperref[app]{Appendix}.

\begin{lem} \label{lem_fd}
Let $ \alpha\in(0, 1/2 ] $ and define
\[
s = s_\nu( x ) = 1-\nu^{-1/2} (-\ln ( x ) )^{1/2} + o(
\nu^{-1/2}), \qquad x \in(0,1/2].
\]
Then, given model D-EX-t$(\zeta_n)$ with $ \lim_{n \to\infty} \zeta _n
= \zeta=1 $, it holds for sufficiently large $\nu$ that $ F_\infty(
\cdot| s_\nu(x) ) $ has (\textup{i}) two CP's for all $ x \in (0,\alpha) $, and
(\textup{ii}) no CP for all $ x \in(\alpha,1/2] $.
\end{lem}

Application of this lemma yields the same limit of the FDR for $
\zeta_n \to1 $ and $ \nu\to\infty$ as in Theorem \ref{thm4_3}; cf. the
discussion in Section~\ref{sec6}.

\begin{thm} \label{F_nu}
Let $ \alpha\in(0, 1/2 ] $. Then, given model D-EX-t$(\zeta_n)$ with\break $
\lim_{n \to\infty} \zeta_n =1 $,
%
%e22 ###
%
\begin{equation}
\lim_{\nu\to\infty} \operatorname{FDR}_{\infty} ( 1 ) = \Phi
\bigl(-\sqrt{-2 \ln(\alpha)} \bigr) . \label{surprise2}
\end{equation}
\end{thm}

\begin{pf}
The result follows by letting $ x \to\alpha$ in Lemma \ref{lem_fd} and
by applying the central limit theorem. Setting $s_\nu= s_\nu( \alpha)
$, we get
\begin{eqnarray*}
\lim_{\nu\to\infty} \operatorname{FDR}_{\infty} ( 1 ) &=&
\lim_{\nu\to\infty} P ( S \leq s_\nu) %\\
%&=&
\\
&=& \lim_{\nu\to\infty} P \biggl( \frac{\nu S^2 - \nu}{\sqrt{2 \nu} }
\leq\frac{\nu s_\nu^2 - \nu}{\sqrt{2 \nu} } \biggr)
\\
&=& \lim_{\nu\to\infty} P \biggl( \frac{\nu S^2 - \nu}{\sqrt{2 \nu} }
\leq-\sqrt{-2 \ln(\alpha)} + o(1) \biggr)
\\
&=& \Phi\bigl(-\sqrt{-2 \ln(\alpha)} \bigr).
\end{eqnarray*}
\upqed
\end{pf}

%f4 ###
%
\begin{figure}[t] % figuur 1

\includegraphics{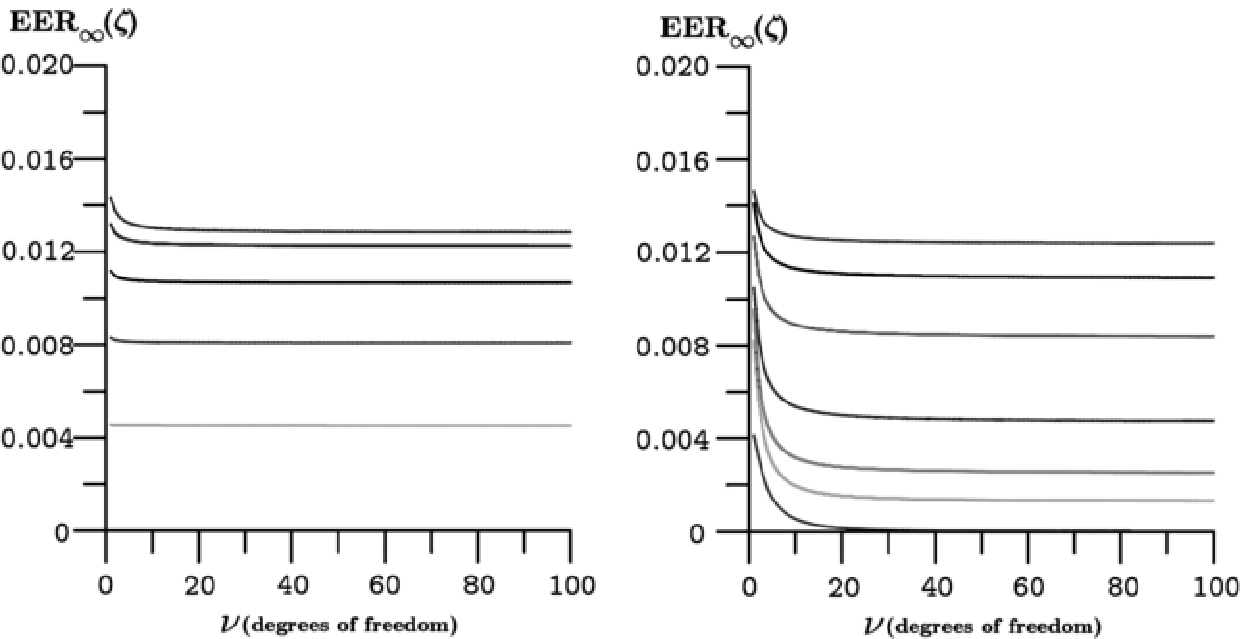}

\caption{$ \operatorname{EER}_\infty(\zeta)$ in the
D-EX-t$(\zeta_n)$ model for $ \alpha=0.05$ and
$\zeta=0.1,~0.2,~0.3,~0.4,~0.5$ (left picture) and
$\zeta=0.6,~0.7,~0.8,~0.9,~0.95,~0.975,~1$ (right picture).}\label{figure_nu_eer}
\end{figure}

\begin{figure}[b]%[H] % figuur 1

\includegraphics{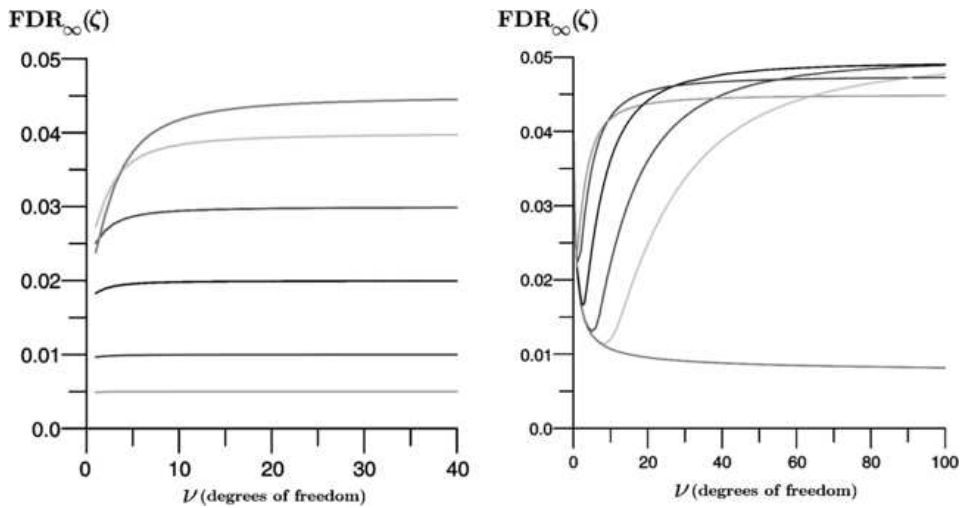}

\caption{$ \operatorname{FDR}_\infty(\zeta)$ in the
D-EX-t$(\zeta_n)$ model for $ \alpha=0.05$ and
$\zeta=0.1, 0.2,0.4,0.6,0.8,0.9$ (left picture) and
$\zeta=0.9,~0.95,~0.99,~0.999,~0.9999,~1$ (right picture).}\label{figure_nu}
\end{figure}

In analogy to Section \ref{sec4}, Figures \ref{figure_nu_eer} and
\ref{figure_nu} display $ \operatorname{EER}_\infty( \zeta) $ and $
\operatorname {FDR}_\infty( \zeta) $, respectively, for various values
of $ \zeta$ and~$ \nu$. It seems that $ \operatorname{EER}_\infty(
\zeta) $ is decreasing in~$ \nu$. For $ \nu\to\infty$, $
\operatorname{EER}_\infty( \zeta) $ again tends to the value $
\alpha(1-\zeta)/(1-\alpha\zeta) $ as expected; see~\cite{r9}.
Note that $ \operatorname{EER}_\infty( \zeta) $ is already close to
this limit if $\nu$ is not too small. As expected, for $
\zeta\approx1/2 $ and $ \nu$ not too small, $
\operatorname{EER}_\infty( \zeta) $ is largest. Except for $ \zeta=1 $,
the FDR tends to the Benjamini--Hochberg bound $\zeta \alpha$ for
increasing degrees of freedom. The limit for $\nu\to0 $ is not clear.
In the latter case the density of the \mbox{$t$-distribution} becomes flatter
and flatter\ and the computation of $ \operatorname{FDR}_\infty( \zeta)
$ becomes extremely difficult. As in the D-EX-N$(\zeta_n)$ model with
$\zeta_n \to1$, $ \operatorname{FDR}_\infty( 1 ) $ reflects the
limiting behavior of the true level of Simes' \cite{r17} global test
for the intersection hypothesis and again, our results show
that it is asymptotically greater than zero for all $\nu> 0$. %\bigskip

%f5 ###
%

%s5 ###
\section{Concluding remarks} \label{sec6}

The investigations in this paper show that the false discovery
proportion FDP$= V_n / [ R_n \vee1 ] $ of the LSU-procedure can be very
volatile in the case of dependent $p$-values, that is, the actual FDP may be
much larger (or smaller) than in the independent case. The same is true
for $ V_n $, $ V_n /n $, $R_n$ and~$R_n/n$. Under mild assumptions, the
e.c.d.f. of the $p$-values converges to a fixed curve under independence
(cf.\ \cite{r9}), which implies convergence of $ V_n /n $ and $ R_n/n $
to fixed values. On the other hand, the shape of the e.c.d.f. of the
$p$-values under exchangeability heavily depends on the (realization of
the) disturbance variable $Z$; cf. Figure \ref{figure_e.c.d.f._Simes}.
In the latter case, the limit distribution of $ V_n /n $ and $ R_n / n
$ typically has positive variance.

It is often assumed that there may be some kind of weak dependence
between test statistics (cf., e.g., \cite{r22}), being close to
independence in some sense. The results in Theorems \ref{thm4_3} and
\ref{F_nu} and the numerical calculations reflected in Figures
\ref{figure_rho} and \ref{figure_nu} suggest that for large $n$ and
$\zeta_n \rightarrow1$ small deviations from independence (small $\rho$
or large $\nu$) may result in a substantially smaller FDR than the
Benjamini--Hochberg bound.
However, simulations for small $ \rho$ and large $ \nu$ show that
$\operatorname{FDR}_{n} ( \zeta_n )$ approaches its limit
$\operatorname{FDR}_{\infty} ( 1 ) $ only for unrealistically large
values of $n$ if \mbox{$ \zeta_n \rightarrow1$}. For example, in the
D-EX-N$(\zeta_n)$ model with $\alpha=0.05$, $ n = 100{,}000$ and $ \rho=
0.1 $, we obtained $ \operatorname{FDR}_{n} ( 1 ) \approx0.0417$ by
simulation. For $\rho=0.01$, we got $ \operatorname{FDR}_{n} ( 1 )
\approx0.05$.
A possible explanation may be that $ \lim_{\rho\rightarrow0^+}
\operatorname{FDR}_n( 1 ) = \alpha$, \mbox{$ \lim_{\nu\rightarrow\infty}
\operatorname{FDR}_n( 1 ) = \alpha$}, hence, the order of limits plays a
severe role. Moreover, for small $ \rho$, it seems that $n$ has to be
very large such that the e.c.d.f. reproduces the shape of $F_\infty$
close to $0$.
For $ \zeta< 1 $, the $\operatorname{FDR}_\infty$ curves in Figures
\ref{figure_rho} and \ref{figure_nu} reflect the FDR for realistically
large $n$ (e.g., $n=1000$) very well. The reason is that the shape
behavior of $F_\infty$ close to $0$ is not that crucial as for $
\zeta=1$.

Example \ref{examp_expo} shows that the FDR under dependence may also
have the same behavior as in the independent case. Therefore, it seems
very difficult to predict what happens with EER, FDR and FDP in models
with more complicated dependence structure, for example, in a multivariate
normal model with arbitrary covariance matrix. In any case, results of
the LSU-procedure, or more generally, of any FDR controlling procedure,
should be interpreted with some care under dependence, taking into
account that the FDR refers to an expectation and that the procedure at
hand may lead to much more false discoveries than expected.

In the models studied in Sections \ref{sec4} and \ref{sec5}, the EER
becomes smallest if $ \vartheta_i \rightarrow0^+ $ for all $ i \in I_1
$ and tends to $ \zeta_n \operatorname{EER}_n ( 1 ) $, where $ I_1 = \{
j \dvtx  K_j \ni\vartheta_j \} $. It is not clear for which parameters
$\vartheta_i$ the FDR becomes smallest. However, if $ \vartheta_i
\rightarrow0^+ $ for all $ i \in I_1 $, $ \zeta_n
\rightarrow\zeta\in(0,1)$, the FDR tends to $ \zeta
\operatorname{FDR}_\infty(1)$.

Finally, with slight modifications of the methods developed in this
paper, one can also treat statistics like $T_i=|X_i-Z|$ or
$T_i=|X_i|/Z$. Somewhat more effort will be necessary if the
disturbance variable $Z$ is two-dimensional as, for example, in $ T_i =
|X_i-Z_1|/Z_2$.

\begin{appendix}

%s6 ###
\section*{Appendix: Proofs}\label{app}

\begin{pf*}{Proof of Theorem \ref{thm2_3}}
The assumptions concerning $ F_\infty$ imply that
$ \lim_{n\to\infty} R_n (z) /n = 0 $ almost surely.
Noting that $ V_n (z) /n \leq R_n (z) /n $ for all $ n \in\mathbb{N}$,
(\ref{vn3}) is obvious.

In order to prove (\ref{vnrn4}), we nest $F_\infty$ between two
c.d.f.'s being linear in a neighborhood of zero. To this end, let $t^*
\in(0, \alpha]$ be fixed, $ B = [0, t^*)$, $ m_\ell(t^*) = \inf_{t \in
B \setminus\{0\}} F_{\infty}(t|z) / t$, $ m_u(t^*) = \sup_{t \in B
\setminus\{0\}} F_{\infty}(t|z) / t $, and ~
\begin{eqnarray*}
F_\ell(t) &=& m_\ell(t^*) t \cdot\mathbf{1}_{B}(t) + F_{\infty
}(t|z) \cdot\mathbf{1}_{B^c}(t),
\\
F_u(t) &=& m_u (t^*) t \cdot\mathbf{1}_{B}(t) + \max\{ m_u (t^*) t^* ,
F_\infty( t | z ) \} \cdot\mathbf{1}_{B^c}(t) .
\end{eqnarray*}
This results in $ F_\ell(t) \leq F_{\infty}(t|z) \leq F_u(t)$ for all $
t \in[0, 1] $.
%We first show (\ref{vnrn4}).
For $n \in\mathbb{N}$,
let the event $A_n(t^*)$ be defined as in Lemma \ref{lem3_4}.
Then
\begin{eqnarray*}
\operatorname{FDR}_n(\zeta_n|z) &= &\mathbb{E}\biggl( \frac
{V_n(z)}{R_n(z) \vee1} \mathbf{1}_{A_n(t^*)} \biggr) + \mathbb
{E}\biggl( \frac{V_n(z)}{R_n(z) \vee1} \mathbf{1}_{A^c_n(t^*)}
\biggr)
\\
&= & \Lambda_n + \lambda_n \mbox{ (say)}.
\end{eqnarray*}
With $r_n = \max\{i \in\mathbb{N}_0\dvtx  i \alpha/ n \leq t^* \}$, we
obtain similarly to the arguments in the proof of Lemma
\ref{FDR_lem1} that
\begin{eqnarray*}
\Lambda_n &= &\mathbb{E}\biggl( \frac{V_n(z)}{R_n(z) \vee1}
\mathbf{1}_{\{R_n(z) \leq r_n\}} \biggr)
\\
&= &n_0 \sum_{j=1}^{r_n} \frac{\mathrm{P}(p_1(z) \leq j\alpha/n)}{j}
\bigl[\mathrm{P}\bigl(D_j^{(1)}(\alpha)\bigr) - \mathrm{P}\bigl(D_{j-1}^{(1)}(\alpha
)\bigr)\bigr].
\end{eqnarray*}
Due to the pointwise order of $ F_\ell$, $ F_\infty$ and $ F_u $, we get
\begin{eqnarray*}
\zeta_n m_{\ell}(t^*) \alpha\mathrm{P}\bigl(D_{r_n}^{(1)}(\alpha)\bigr) & \leq&
\Lambda_n \leq \zeta_n m_u(t^*) \alpha\mathrm{P}\bigl(D_{r_n}^{(1)}(\alpha)\bigr),
\\
\zeta_n m_{\ell}(t^*) \alpha\mathrm{P}\bigl(D_{r_n}^{(1)}(\alpha)\bigr) +
\lambda_n & \leq& \operatorname{FDR}_n(\zeta_n|z) \leq \zeta_n
m_u(t^*) \alpha\mathrm{P}\bigl(D_{r_n}^{(1)}(\alpha)\bigr) + \lambda_n.
\end{eqnarray*}
Since $\zeta_n \to1$, $\mathrm{P}(D_{r_n}^{(1)}(\alpha)) \to1$ and
$\mathrm{P}(A_n(t^*)) \to1$ for $n \to\infty$, we obtain \mbox{$ \lambda_n
\to0$} and $m_\ell(t^*) \alpha\leq\lim\inf_{n \to\infty} \operatorname
{FDR}_n(\zeta_n|z) \leq\lim\sup_{n \to\infty} \operatorname
{FDR}_n(\zeta_n|z)$ $ \leq\break m_u(t^*) \alpha$. The assertion now follows
by noticing that $\lim_{t^* \to0^+} m_\ell (t^*) =\break \lim_{t^* \to0^+}
m_u(t^*) = \gamma(z)$.
\end{pf*}

\begin{pf*}{Proof of Theorem \ref{thm4_1}}
Denote the p.d.f. corresponding to $ G_{\zeta,1} $ by $ g_{\zeta,1} $
and notice that $C_{\zeta, 1} = (0, t_1 / \alpha- (1- \zeta)) \cup (t_2
/ \alpha- (1- \zeta), \zeta)$. From Theorem~\ref{thm3_2}, we get
\[
\operatorname{EER}_\infty( \zeta)
= \int_0^{t_1 / \alpha- (1- \zeta)} u g_{\zeta,1}(u)\, du + \int_{t_2
/ \alpha- (1- \zeta)}^\zeta u g_{\zeta,1}(u) \,du.
\]
Since $x_0(t_1 | \zeta) = x_0(t_2 | \zeta)$ and $\lim_{t \uparrow
\alpha} x_0(t|\zeta) = -\infty$, we get
\begin{eqnarray*}
\operatorname{EER}_\infty( \zeta) &= & \bigl(t_1 / \alpha- (1- \zeta)\bigr)
\bigl(1-\Phi(x_0(t_1|\zeta))\bigr) +\zeta
\\
&&{} - \bigl(t_2 / \alpha- (1- \zeta)\bigr) \bigl(1-\Phi(x_0(t_1|\zeta))\bigr) - \bigl(t_1
/ \alpha- (1- \zeta)\bigr) - \zeta
\\
&&{} + t_2 / \alpha- (1- \zeta) + \int_0^{t_1 / \alpha- (1- \zeta
)} \Phi\bigl(x_0\bigl(\alpha(u+1-\zeta)|\zeta\bigr)\bigr)\, du
\\
&&{} +\int_{t_2 / \alpha- (1- \zeta)}^\zeta\Phi\bigl(x_0\bigl(\alpha(u+1-\zeta
)|\zeta\bigr)\bigr) \,du
\\
& = & \frac{t_2 - t_1}{ \alpha} \Phi( x_0 ( t_1 | \zeta) )
\\
&&{} + \int_{1-\zeta}^{t_1/\alpha} \Phi( x_0 ( \alpha t | \zeta) ) \,dt
+ \int_{t_2 / \alpha}^1 \Phi( x_0 ( \alpha t | \zeta) )\,dt.
\end{eqnarray*}
In order to compute $\operatorname{FDR}_\infty(\zeta)$, note that, for
$z \in( 0 , z_1 ) \cup( z_2 , \zeta)$,
\[
G_{\zeta,2} ( z )
= 1 - \Phi\biggl( x_0 \biggl( \frac{\alpha( 1 - \zeta)}{1 - z } \Big| \zeta\biggr) \biggr) ,
\]
where $ z_i = 1- \alpha( 1 - \zeta) / t_ i $ , $ i = 1 , 2 $. In view
of $\lim_{t \downarrow\alpha(1 - \zeta)} x_0(t|\zeta) = \infty$, it is
$ G_{\zeta,2} ( z_1 ) = G_{\zeta,2} ( z_2 ) $, $G_{\zeta,2} ( 0 ) = 0 $
and $ G_{\zeta,2} (\zeta)=1$. Denoting the corresponding p.d.f. of
$G_{\zeta,2}$ by $ g_{\zeta,2} $, we obtain
\begin{eqnarray*}
\operatorname{FDR}_{\infty} ( \zeta)
& = & \int_{ 0 }^{ z_1 } z g_{\zeta,2} ( z ) \,dz + \int_{ z_2 }^{ \zeta} z
g_{\zeta,2} ( z ) \,dz
\\
& = & z_1 G_{\zeta,2} ( z_1) + \zeta G_{\zeta,2} ( \zeta) - z_2 G_{\zeta
,2} ( z_2)
\\
&&{} - \int_{ 0 }^{ z_1 } G_{\zeta,2} ( z ) \,d z - \int_{ z_2 }^{
\zeta} G_{\zeta,2} ( z ) \,d z
\\
& = & ( z_2 - z_1 ) \Phi\biggl( x_0 \biggl( \frac{\alpha( 1 - \zeta)}{1 - z_1 } \Big|
\zeta\biggr) \biggr) \\
&&{} + \int_{ 0 }^{ z_1 } \Phi\biggl( x_0 \biggl( \frac{\alpha( 1 - \zeta)}{1 - z }
\Big| \zeta\biggr) \biggr)\, d z
\\
&&{}+ \int_{ z_2 }^{ \zeta} \Phi\biggl( x_0 \biggl( \frac{\alpha( 1 - \zeta)}{1 -
z } \Big| \zeta\biggr) \biggr) \,d z .
\end{eqnarray*}
\upqed
\end{pf*}

\begin{pf*}{Proof of Theorem \ref{thm4_3}}
For any $ \rho\in(0,1) $, there exists a unique solution
$(u,x_0)=(u_\rho,x_{0,\rho})$ (say) of (\ref{BP_cond_1}) and
(\ref{BP_cond_2}).
%order to calculate the limiting value for $x_0$ leading to a TP in
%case of $\rho\to0^+$
%for $ \zeta= 1 $, we utilize that the smaller solution of (
In view of (\ref{BP_cond_3}) and the shape of $F_\infty$, $(u_\rho,x_{0,\rho})$ satisfies
%
%e23 ###
%
\begin{equation} \label{ideales_u}
%u^*( x_0 ) = - x_0 / \sqrt{ \rho} - \sqrt{ \rhoq/ \rho}
u_\rho= - x_{0,\rho} / \sqrt{ \rho} - \sqrt{ \overline{\rho}/
\rho}
\sqrt{ x_{0,\rho}^2 - 2 \ln ( \sqrt{ \overline{\rho}}/ \alpha )}.\vadjust{\goodbreak}
\end{equation}
Notice that $ \alpha\in( 0 , 1/2 ] $ implies $ u_\rho> 0 $ and
therefore, $ x_{0 , \rho} < 0 $.
Now, for $\delta\in(0, \alpha)$, we consider $x_0 = x_0(\delta) =
-\sqrt{- 2 \ln(\delta)} < -\sqrt{- 2 \ln(\alpha)} = x_0(\alpha)$ in
order to bound $x_{0,\rho}$ from below for $ \rho\to0^+ $. Since
$u_\rho$ has to be a real number, we obtain from (\ref{ideales_u}) that
$\limsup_{\rho\to0^+} x_{0,\rho} \leq x_0(\alpha)$.
%only $\delta\leq\alpha$ can lead to a solution.
Defining
\[
u = u(\rho, \delta) = \frac{-x_0(\delta)}{\sqrt{\rho}} \quad\mbox{and}\quad
w = w(\rho, \delta) = \frac{u(\rho, \delta)}{\sqrt
{\overline{\rho}}} + \sqrt{\frac{\rho}{\overline{\rho}}}
x_0(\delta),
\]
we get from (\ref{distanz}) that $ d(u|x_0,1) = \Phi(-w) - \Phi(-u) /
\alpha. $ Hence, $ d(u|x_0,1) > 0 $ iff \mbox{$ \Phi(-u) / \Phi(-w) <
\alpha$}. Employing the asymptotic relationship $
\Phi(-x) / \varphi(x) \sim1 / x\ (x \to\infty)$ for Mills' ratio, we get
\[
\frac{\Phi(-u)}{\Phi(-w)} \sim\frac{w}{u} \frac{\varphi
(u)}{\varphi(w)}
= \frac{w}{u} \exp\bigl((w^2 - u^2) / 2\bigr).
\]
Since $\exp((w(\rho, \delta)^2 - u( \rho, \delta)^2 ) / 2) = \delta<
\alpha$ independent of $\rho$ and $\lim_{\rho\to0^+} w(\rho,\break \delta) /
u(\rho, \delta) = 1$, we conclude that $\lim_{\rho\to0^+} x_{0,\rho} =
x_0(\alpha) =- \sqrt{ - 2 \ln (\alpha) }$. This finally implies
(\ref{surprise1}) and completes the proof.
\end{pf*}

\begin{pf*}{Proof of Lemma \ref{lem_fd}}
For $ s^2 < (\nu+1)/\nu$, the unique point of inflection of
$F_\infty(\cdot|s)$ on $ (0,1/2)$ is given by $ t^*(\nu|s) = F_{t_\nu}
(-a(s,\nu)) $ with $ a(s,\nu)=\break((\nu+1)/s^2-\nu)^{1/2}$. Hence, it
suffices to show that
\[
F_\infty\bigl( t^*(\nu|s_\nu(x)) | s_\nu( x ) \bigr) > t^*(\nu|s_\nu( x )
)/ \alpha \qquad\mbox{for } x\in(0,\alpha)
\]
for sufficiently large $\nu$ and that the derivative of $ F_\infty(
\cdot| s_\nu(x)) $ in $ t = t^*(\nu|s_\nu( x ) ) $ is less than $ 1 /
\alpha$ for all $ x \in( \alpha, 1/2 ] $ for sufficiently large $\nu$.
Therefore, the assertion follows if
%e25 ###
%e24 ###
%
\begin{eqnarray} \label{yeah1}
\lim_{ \nu\to\infty} \frac{F_{t_\nu} ( - a( s_\nu( x ) , \nu) )
}{ \Phi( - s_\nu( x ) a( s_\nu( x ) , \nu) )}
& < & \alpha \qquad\mbox{for } x\in(0,\alpha),
\\
\lim_{ \nu\to\infty}
\frac{ f_{t_\nu}  ( a( s_\nu( x ) , \nu) ) }{ s_\nu( x ) \varphi
 ( s_\nu( x ) a( s_\nu( x ) , \nu) ) }
& > & \alpha \qquad\mbox{for } x\in(\alpha,1/2] . \label{yeah2}
\end{eqnarray}
For $ x_\nu\in( 0 , \infty) $ with $ \lim_{ \nu\to\infty} x_\nu^4 /
\nu= \beta\in[0,\infty]$, it is shown in \cite{r4fd} that
\[
\lim_{\nu\to\infty} \frac{ f_{t_\nu} ( x_\nu) }{ \varphi( x_\nu )} =
\lim_{\nu\to\infty} \frac{ F_{t_\nu} ( - x_\nu ) }{ \Phi( - x_\nu)} =
\exp(\beta/ 4).
\]
Note that for $ u \to\infty$ and $ s \to1 $, it holds that (Mills' ratio)
\[
\frac{ F_{t_\nu} ( - u ) }{ \Phi( -s u )} \sim\frac{ F_{t_\nu} ( - u )
}{ \Phi( - u )} \frac{ \varphi( u ) }{ \varphi( s u )}.
\]
We\vadjust{\goodbreak} easily get $ \lim_{\nu\to\infty} a(s_\nu(x),\nu)^4/\nu =
\lim_{\nu\to\infty} a(s_\nu(x),\nu)^2(1-s_\nu(x)^2) = - 4 \ln( x ) $.
As a consequence, (\ref{yeah1}) follows by noting that
\begin{eqnarray*}
&&\lim_{ \nu\to\infty} \frac{F_{t_\nu} ( - a( s_\nu( x ) , \nu) ) }{
\Phi( - s_\nu( x )a( s_\nu( x ) , \nu) )}
\\
&&\qquad = \lim_{ \nu\to\infty}
\biggl[
\frac{F_{t_\nu} ( - a( s_\nu( x ) , \nu) ) }{ \Phi( - a( s_\nu( x
) , \nu) )}
\frac{ \varphi( a( s_\nu( x ) , \nu) ) }{ \varphi( s_\nu( x )
a( s_\nu( x ) , \nu) )}
\biggr]
\\
&&\qquad = \exp(- 4 \ln( x ) / 4 )
\lim_{ \nu\to\infty}\exp\biggl( - \frac{1}{2} a( s_\nu( x ) , \nu)^2 \bigl(1- s_\nu( x )^2\bigr)
\biggr)
\\
&&\qquad= \frac{1}{x}\exp(2 \ln(x)) = x.
\end{eqnarray*}
An analogous calculation yields (\ref{yeah2}). Hence, Lemma
\ref{lem_fd} is proved.
\end{pf*}

\end{appendix}

\section*{Acknowledgments}
The authors are grateful to a referee and an Associate Editor for their
constructive and valuable suggestions and their quick replies. Thanks
are also due to the editor M. L. Eaton for his expeditious handling
of the manuscript.

\printaddresses

\end{document}